\documentclass[a4paper,twoside,11pt,english,fleqn]{article}

\usepackage[all]{xy}
\usepackage[intlimits]{amsmath}
\usepackage{amssymb,amsthm,graphicx}
\usepackage[dvips,dvipdf,pdftex]{epsfig}
\usepackage{babel}
\usepackage[T1]{fontenc}
\usepackage{fancyhdr,stmaryrd,color,times,euler,euscript,eufrak}

\DeclareGraphicsExtensions{.png}
\DeclareMathAlphabet{\mathscr}{T1}{pzc}{m}{it}

  \newtheorem{thm}{Theorem}[section]
  \newtheorem{cor}[thm]{Corollary}
  \newtheorem{lem}[thm]{Lemma}
  \newtheorem{prop}[thm]{Proposition}
\theoremstyle{definition}
  \newtheorem{defn}[thm]{Definition}
  \newtheorem{ex}[thm]{Example}
  \newtheorem{rem}[thm]{Remark}
  \newtheorem{ntt}[thm]{Notation}
  \newtheorem*{dem}{Proof}

\newcommand{\Cr}{\EuScript{C}} 
\newcommand{\Gr}{\EuScript{G}} 
\newcommand{\Mr}{\EuScript{M}} 
\newcommand{\Nr}{\EuScript{N}} 
\newcommand{\Nb}{\mathbb{N}} 
\newcommand{\Sk}{\EuFrak{S}}
\renewcommand{\phi}{\varphi}
\renewcommand{\epsilon}{\varepsilon}
\newcommand{\findem}{\hfill $\diamondsuit$}

\newcommand{\fl}{\rightarrow}

\newcommand{\mfl}{\twoheadrightarrow}

\newcommand{\red}[1]{\rightarrow\!\!_{{\scriptscriptstyle #1}}}
\newcommand{\mred}[1]{\twoheadrightarrow\!\!_{{\scriptscriptstyle #1}}}
\newcommand{\equi}[1]{\equiv\!\!_{{\scriptscriptstyle #1}}}

\newcommand{\tens}{\otimes}
\newcommand{\mon}[1]{\langle #1 \rangle}
\newcommand{\ens}[1]{\{#1\}}
\newcommand{\dens}[2]{\{#1,\dots,#2\}}
\newcommand{\ol}[1]{\overline{#1}}

\newcommand{\et}{\quad\text{and}\quad}

\DeclareMathOperator{\id}{id}

\newcommand{\emptysectionmark}[1]
	{\markboth{\textbf{#1}}{\textbf{#1}}}

\fancypagestyle{plain}
{
  \fancyhf{}
	\fancyfoot[LC,RC]{}
  \fancyfoot[LE,RO]{\textbf{\thepage}}

}

\begin{document}
\thispagestyle{plain}

\hfill {\large \textbf{1st July 2005 - Modified: 6th January 2006}} 

\vspace{3mm}
\hrule height 1.5pt

\vspace{4mm}
{\LARGE \textbf{TWO POLYGRAPHIC PRESENTATIONS}}

\vspace{1mm}
{\LARGE \textbf{OF PETRI NETS}}

\vspace{3mm}
\indent{\LARGE \textbf{Yves Guiraud}\footnote{Institut de mathématiques de Luminy, Marseille, France - \texttt{guiraud@iml.univ-mrs.fr}} } 

\vspace{3mm}
\hrule height 1.5pt

\vspace{8mm}

\begin{minipage}{150mm}
\textbf{Abstract:} This document gives an algebraic and two polygraphic translations of Petri nets, all three providing an easier way to describe reductions and to identify some of them. The first one sees places as generators of a commutative monoid and transitions as rewriting rules on it: this setting is totally equivalent to Petri nets, but lacks any graphical intuition. The second one considers places as $1$-dimensional cells and transitions as $2$-dimensional ones: this translation recovers a graphical meaning but raises many difficulties since it uses explicit permutations. Finally, the third translation sees places as degenerated $2$-dimensional cells and transitions as $3$-dimensional ones: this is a setting equivalent to Petri nets, equipped with a graphical interpretation.
\end{minipage}

\section*{Outline}
\emptysectionmark{Outline}

In this document, we study Petri nets in order to give two possible polygraphic presentations for them. This work follows Albert Burroni's intuitions: many computer science and proof theory objects have natural translations into polygraphs. These are topology-flavoured objects consisting of collections of directed cells of various dimensions, equipped with a rich algebraic structure.

In section~\ref{sec:reseaux_de_petri}, we recall some basic facts about Petri nets, describe their representations and associate them reduction graphs, equipped with a relation that identifies paths that intuitively represent the same sequence of operations.

In section~\ref{sec:mots_commutatifs}, we recall a known algebraic account of Petri nets: they correspond to commutative word rewriting systems (or presentations of commutative monoids) and both objects generate the same reduction graph. Furthermore, in the latter, reductions have a name, which makes easier the definition of a relation between similar paths. We prove a new result concerning stating that this relation is the same as the one defined for Petri nets. All these facts are detailed in theorem~\ref{th:comm}.

In section~\ref{sec:2d}, we craft a $2$-dimensional object, a $2$-polygraph, in which reductions of a Petri net can be translated. This result is due to Albert Burroni and is formulated as theorem~\ref{th:2d}. We go beyond and study the links between the relation on Petri nets paths and two relations on $2$-arrows of the $2$-polygraph: the first one corresponds to the relation on the Petri net, while the second one tries to solve the difficulties raised by the presence of explicit permutations in the $2$-polygraph. The study of these properties is only started here: much more work will be necessary to totally solve the encountered problems. 

Finally, in section~\ref{sec:3d}, we give a new, more natural polygraphic way to faithfully describe Petri nets. We prove that they correspond to $3$-polygraphs with one cell in dimension~$0$ and no cell in dimension~$1$. Furthermore, both objects generate the same reduction graph, with the same equivalence relation on paths. This is the main result, theorem~\ref{th:3d}.

\section{Basic notions on Petri nets}
\label{sec:reseaux_de_petri}

\noindent This section briefly recalls the basic notions about Petri nets: the definitions of a net, of its markings and the usual associated graphical representations. It should be noted that there exist many possible definitions of Petri nets, but a simple one has been chosen for this study. More of them can be found in [Murata 1989] for example. 

\begin{defn}\label{def:reseau_de_petri}
A \emph{Petri net} is a quadruple $N=(X,T,w,w')$ made of two finite sets, $X$ and $T$, and two maps, $w:X\times T\fl\Nb$ and $w':T\times X\fl\Nb$. The elements of $X$ and $T$ are respectively called \emph{places} and \emph{transitions}, while the maps $w$ and $w'$ are the \emph{weights}. Beside this set-theoretic definition, Petri nets are usually encountered as graphical objects. A decorated graph is associated to a given net $N=(X,T,w,w')$ as follows: 
\begin{enumerate}
\item[0.] Its objects are the places and the transitions. Places are pictured as circles, while transitions are represented by double bars.
\item[1.] If $x$ is a place and $\alpha$ a transition, there is an arrow from $x$ to $\alpha$ whenever $w(x,\alpha)>0$ and one from $\alpha$ to $x$ whenever $w'(\alpha,x)>0$. Such arrows are decorated with the corresponding weight, either $w(x,\alpha)$ or $w'(\alpha,x)$.
\end{enumerate}
\end{defn}

\begin{ex}\label{ex:reseau_de_petri}
Let us condider the Petri net $N=(X,T,w,w')$ where $X=\ens{x,y,z}$, $T=\ens{\alpha,\beta}$ and the non-zero values of $w$ and $w'$ are given by:
$$
w(x,\alpha) \:=\: 1, \quad w(y,\beta) \:=\: 2, \quad w'(\alpha,y) \:=\: w'(\alpha,z) \:=\: w'(\beta,z)\: =\: 1.
$$

\noindent Following the given graph construction recipe, this representation is built for $N$:
\begin{center}
\begin{picture}(0,0)%
\includegraphics{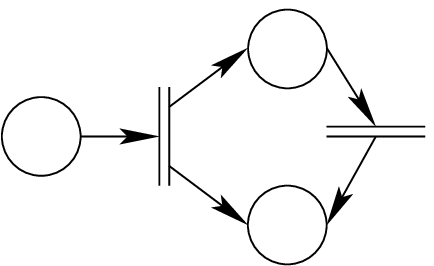}%
\end{picture}%
\setlength{\unitlength}{4144sp}%
\begingroup\makeatletter\ifx\SetFigFont\undefined%
\gdef\SetFigFont#1#2#3#4#5{%
  \reset@font\fontsize{#1}{#2pt}%
  \fontfamily{#3}\fontseries{#4}\fontshape{#5}%
  \selectfont}%
\fi\endgroup%
\begin{picture}(2078,1179)(173,-518)
\put(1891,-241){\makebox(0,0)[b]{\smash{{\SetFigFont{10}{12.0}{\rmdefault}{\mddefault}{\updefault}{\color[rgb]{0,0,0}$1$}%
}}}}
\put(361, 29){\makebox(0,0)[b]{\smash{{\SetFigFont{12}{14.4}{\rmdefault}{\mddefault}{\updefault}{\color[rgb]{0,0,0}$x$}%
}}}}
\put(1486,434){\makebox(0,0)[b]{\smash{{\SetFigFont{12}{14.4}{\rmdefault}{\mddefault}{\updefault}{\color[rgb]{0,0,0}$y$}%
}}}}
\put(1486,-376){\makebox(0,0)[b]{\smash{{\SetFigFont{12}{14.4}{\rmdefault}{\mddefault}{\updefault}{\color[rgb]{0,0,0}$z$}%
}}}}
\put(901,-331){\makebox(0,0)[b]{\smash{{\SetFigFont{12}{14.4}{\rmdefault}{\mddefault}{\updefault}{\color[rgb]{0,0,0}$\alpha$}%
}}}}
\put(2251, 29){\makebox(0,0)[b]{\smash{{\SetFigFont{12}{14.4}{\rmdefault}{\mddefault}{\updefault}{\color[rgb]{0,0,0}$\beta$}%
}}}}
\put(721,119){\makebox(0,0)[b]{\smash{{\SetFigFont{10}{12.0}{\rmdefault}{\mddefault}{\updefault}{\color[rgb]{0,0,0}$1$}%
}}}}
\put(1081,389){\makebox(0,0)[b]{\smash{{\SetFigFont{10}{12.0}{\rmdefault}{\mddefault}{\updefault}{\color[rgb]{0,0,0}$1$}%
}}}}
\put(1171,-151){\makebox(0,0)[b]{\smash{{\SetFigFont{10}{12.0}{\rmdefault}{\mddefault}{\updefault}{\color[rgb]{0,0,0}$1$}%
}}}}
\put(1846,344){\makebox(0,0)[b]{\smash{{\SetFigFont{10}{12.0}{\rmdefault}{\mddefault}{\updefault}{\color[rgb]{0,0,0}$2$}%
}}}}
\end{picture}%
\end{center}
\end{ex}

\noindent So far, only the hardware part of a Petri net has been represented. On top of this one, the states of the automaton are described:

\begin{defn}\label{def:marquage}
Let $N=(X,T,w,w')$ be a Petri net. A \emph{marking} of $N$ is a map from the set $X$ of places to the set $\Nb$ of natural numbers. The set of all markings of $N$ is denoted by $\Mr(N)$. A given marking $\mu:X\fl\Nb$ on a Petri net $N=(X,T,w,w')$ is represented as an extra decoration on the corresponding graph: inside each place $x$, one puts $\mu(x)$ token(s), pictured as black dots.
\end{defn}

\begin{ex}\label{ex:marquage}
With the same Petri net as in example \ref{ex:reseau_de_petri}, the marking $\mu$ defined by $\mu(x)=\mu(y)=2$ and $\mu(z)=0$ is represented as follows (thereafter, the weights equal to $1$ are removed, together with places labels, in order to make the representations clearer):
\begin{center}
\begin{picture}(0,0)%
\includegraphics{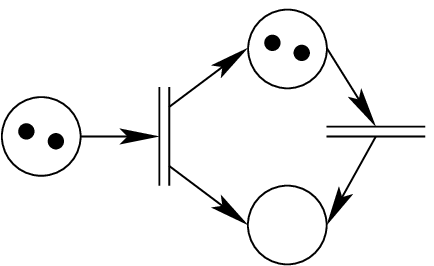}%
\end{picture}%
\setlength{\unitlength}{4144sp}%
\begingroup\makeatletter\ifx\SetFigFont\undefined%
\gdef\SetFigFont#1#2#3#4#5{%
  \reset@font\fontsize{#1}{#2pt}%
  \fontfamily{#3}\fontseries{#4}\fontshape{#5}%
  \selectfont}%
\fi\endgroup%
\begin{picture}(2078,1179)(173,-518)
\put(1846,344){\makebox(0,0)[b]{\smash{{\SetFigFont{10}{12.0}{\rmdefault}{\mddefault}{\updefault}{\color[rgb]{0,0,0}$2$}%
}}}}
\put(901,-331){\makebox(0,0)[b]{\smash{{\SetFigFont{12}{14.4}{\rmdefault}{\mddefault}{\updefault}{\color[rgb]{0,0,0}$\alpha$}%
}}}}
\put(2251, 29){\makebox(0,0)[b]{\smash{{\SetFigFont{12}{14.4}{\rmdefault}{\mddefault}{\updefault}{\color[rgb]{0,0,0}$\beta$}%
}}}}
\end{picture}%
\end{center}
\end{ex}

\noindent Now, the whole static part of Petri nets has been introduced. Their evolutions are described as follows:

\begin{defn}\label{def:transition}
Let $N=(X,T,w,w')$ be a Petri net and let $\alpha$ be a transition in $T$. The \emph{reduction relation associated to $\alpha$} is the binary relation $\red{\alpha}$ on markings of $N$, defined by $\mu\red{\alpha}\nu$ if, for every place $x$ in~$X$, both following conditions hold:
$$
\begin{cases}
\mu(x)\geq w(x,\alpha), \\
\nu(x)=\mu(x)-w(x,\alpha)+w'(\alpha,x).
\end{cases}
$$

\noindent The union of all the relations $\red{\alpha}$, for all the transitions $\alpha$, is denoted by $\red{T}$. The reflexive and transitive closure of $\red{T}$ is denoted by $\mred{T}$ and called the \emph{reachability relation}. 

The relation $\red{\alpha}$ associated to a transition $\alpha$ has a graphical interpretation. The first condition checks if the marking $\mu$ has at least $w(x,\alpha)$ tokens in each place $x$. In that case, the second condition tells that~$\nu$ is entirely determined this way: in each place $x$, $w(x,\alpha)$ tokens are removed, then $w'(\alpha,x)$ tokens are added.
\end{defn}

\begin{ex}\label{ex:graphe_de_reduction}
Let $N$ be the Petri net of example \ref{ex:reseau_de_petri} and $\mu$ the marking of example \ref{ex:marquage}. The graph pictured thereafter displays all the markings of $N$ that can be reached from $\mu$.
\begin{center}
\begin{picture}(0,0)%
\includegraphics{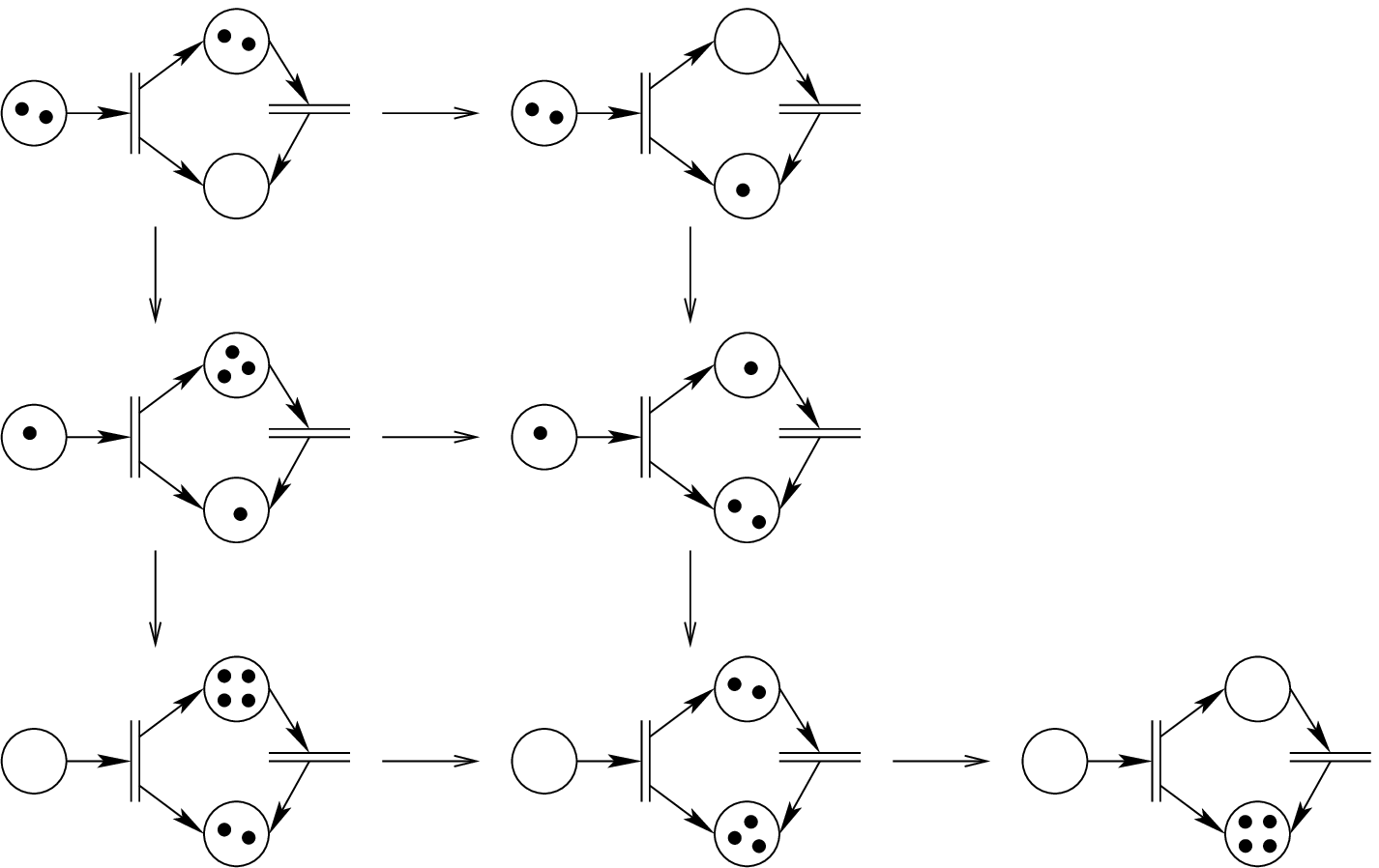}%
\end{picture}%
\setlength{\unitlength}{4144sp}%
\begingroup\makeatletter\ifx\SetFigFont\undefined%
\gdef\SetFigFont#1#2#3#4#5{%
  \reset@font\fontsize{#1}{#2pt}%
  \fontfamily{#3}\fontseries{#4}\fontshape{#5}%
  \selectfont}%
\fi\endgroup%
\begin{picture}(6494,4071)(212,-3717)
\put(6694,-3365){\makebox(0,0)[b]{\smash{{\SetFigFont{9}{10.8}{\rmdefault}{\mddefault}{\updefault}{\color[rgb]{0,0,0}$\beta$}%
}}}}
\put(1675, 83){\makebox(0,0)[b]{\smash{{\SetFigFont{9}{10.8}{\rmdefault}{\mddefault}{\updefault}{\color[rgb]{0,0,0}$2$}%
}}}}
\put(4089, 83){\makebox(0,0)[b]{\smash{{\SetFigFont{9}{10.8}{\rmdefault}{\mddefault}{\updefault}{\color[rgb]{0,0,0}$2$}%
}}}}
\put(1675,-1450){\makebox(0,0)[b]{\smash{{\SetFigFont{9}{10.8}{\rmdefault}{\mddefault}{\updefault}{\color[rgb]{0,0,0}$2$}%
}}}}
\put(4089,-1450){\makebox(0,0)[b]{\smash{{\SetFigFont{9}{10.8}{\rmdefault}{\mddefault}{\updefault}{\color[rgb]{0,0,0}$2$}%
}}}}
\put(1675,-2982){\makebox(0,0)[b]{\smash{{\SetFigFont{9}{10.8}{\rmdefault}{\mddefault}{\updefault}{\color[rgb]{0,0,0}$2$}%
}}}}
\put(4089,-2982){\makebox(0,0)[b]{\smash{{\SetFigFont{9}{10.8}{\rmdefault}{\mddefault}{\updefault}{\color[rgb]{0,0,0}$2$}%
}}}}
\put(6503,-2982){\makebox(0,0)[b]{\smash{{\SetFigFont{9}{10.8}{\rmdefault}{\mddefault}{\updefault}{\color[rgb]{0,0,0}$2$}%
}}}}
\put(832,-454){\makebox(0,0)[b]{\smash{{\SetFigFont{9}{10.8}{\rmdefault}{\mddefault}{\updefault}{\color[rgb]{0,0,0}$\alpha$}%
}}}}
\put(832,-1986){\makebox(0,0)[b]{\smash{{\SetFigFont{9}{10.8}{\rmdefault}{\mddefault}{\updefault}{\color[rgb]{0,0,0}$\alpha$}%
}}}}
\put(832,-3518){\makebox(0,0)[b]{\smash{{\SetFigFont{9}{10.8}{\rmdefault}{\mddefault}{\updefault}{\color[rgb]{0,0,0}$\alpha$}%
}}}}
\put(3246,-454){\makebox(0,0)[b]{\smash{{\SetFigFont{9}{10.8}{\rmdefault}{\mddefault}{\updefault}{\color[rgb]{0,0,0}$\alpha$}%
}}}}
\put(3246,-1986){\makebox(0,0)[b]{\smash{{\SetFigFont{9}{10.8}{\rmdefault}{\mddefault}{\updefault}{\color[rgb]{0,0,0}$\alpha$}%
}}}}
\put(3246,-3518){\makebox(0,0)[b]{\smash{{\SetFigFont{9}{10.8}{\rmdefault}{\mddefault}{\updefault}{\color[rgb]{0,0,0}$\alpha$}%
}}}}
\put(5660,-3518){\makebox(0,0)[b]{\smash{{\SetFigFont{9}{10.8}{\rmdefault}{\mddefault}{\updefault}{\color[rgb]{0,0,0}$\alpha$}%
}}}}
\put(3361,-952){\makebox(0,0)[b]{\smash{{\SetFigFont{9}{10.8}{\rmdefault}{\mddefault}{\updefault}{\color[rgb]{0,0,0}$\alpha$}%
}}}}
\put(871,-2484){\makebox(0,0)[b]{\smash{{\SetFigFont{9}{10.8}{\rmdefault}{\mddefault}{\updefault}{\color[rgb]{0,0,0}$\alpha$}%
}}}}
\put(3361,-2484){\makebox(0,0)[b]{\smash{{\SetFigFont{9}{10.8}{\rmdefault}{\mddefault}{\updefault}{\color[rgb]{0,0,0}$\alpha$}%
}}}}
\put(832,-952){\makebox(0,0)[b]{\smash{{\SetFigFont{9}{10.8}{\rmdefault}{\mddefault}{\updefault}{\color[rgb]{0,0,0}$\alpha$}%
}}}}
\put(2250,-71){\makebox(0,0)[b]{\smash{{\SetFigFont{9}{10.8}{\rmdefault}{\mddefault}{\updefault}{\color[rgb]{0,0,0}$\beta$}%
}}}}
\put(2250,-1603){\makebox(0,0)[b]{\smash{{\SetFigFont{9}{10.8}{\rmdefault}{\mddefault}{\updefault}{\color[rgb]{0,0,0}$\beta$}%
}}}}
\put(2250,-3135){\makebox(0,0)[b]{\smash{{\SetFigFont{9}{10.8}{\rmdefault}{\mddefault}{\updefault}{\color[rgb]{0,0,0}$\beta$}%
}}}}
\put(4664,-3135){\makebox(0,0)[b]{\smash{{\SetFigFont{9}{10.8}{\rmdefault}{\mddefault}{\updefault}{\color[rgb]{0,0,0}$\beta$}%
}}}}
\put(1867,-300){\makebox(0,0)[b]{\smash{{\SetFigFont{9}{10.8}{\rmdefault}{\mddefault}{\updefault}{\color[rgb]{0,0,0}$\beta$}%
}}}}
\put(4281,-300){\makebox(0,0)[b]{\smash{{\SetFigFont{9}{10.8}{\rmdefault}{\mddefault}{\updefault}{\color[rgb]{0,0,0}$\beta$}%
}}}}
\put(4281,-1833){\makebox(0,0)[b]{\smash{{\SetFigFont{9}{10.8}{\rmdefault}{\mddefault}{\updefault}{\color[rgb]{0,0,0}$\beta$}%
}}}}
\put(1867,-1833){\makebox(0,0)[b]{\smash{{\SetFigFont{9}{10.8}{\rmdefault}{\mddefault}{\updefault}{\color[rgb]{0,0,0}$\beta$}%
}}}}
\put(1867,-3365){\makebox(0,0)[b]{\smash{{\SetFigFont{9}{10.8}{\rmdefault}{\mddefault}{\updefault}{\color[rgb]{0,0,0}$\beta$}%
}}}}
\put(4281,-3365){\makebox(0,0)[b]{\smash{{\SetFigFont{9}{10.8}{\rmdefault}{\mddefault}{\updefault}{\color[rgb]{0,0,0}$\beta$}%
}}}}
\end{picture}%
\end{center}
\end{ex}

\vfill\pagebreak
\noindent In order to compare Petri nets with the rewriting-flavoured objects to be introduced in the next three sections, the notion of reduction graph appearing in example \ref{ex:graphe_de_reduction} is formalized:

\begin{defn}\label{def:graphe}
Let $N=(X,T,w,w')$ be a Petri net. Its \emph{associated reduction graph} is the graph $G(N)$ defined by:
\begin{enumerate}
\item[0.] The set of objects of $G(N)$ is the set $\Mr(N)$ of markings of $N$.
\item[1.] In $G(N)$, there is an arrow from a marking $\mu$ to a marking $\nu$ for each transition $\alpha$ such that $\mu\red{\alpha}\nu$.
\end{enumerate}
\end{defn}

\noindent In example \ref{ex:graphe_de_reduction}, we have pictured a subgraph of the reduction graph $G(N)$, where $N$ is the Petri net of example \ref{ex:reseau_de_petri}. Let us consider the top-most square. We can see that the two vertical arrows, both labelled by $\alpha$ are "intuitively" the same reduction: indeed, they consume the same tokens and produce the same ones. This is also the case for the two vertical $\beta$-labelled arrows. Furthermore, the horizontal and vertical reductions apply on different tokens: there should be some relation between the two sequences $\alpha$-then-$\beta$ and $\beta$-then-$\alpha$. Let us define a congruence relation on such reduction paths:

\begin{ntt}
Let $N=(X,T,w,w')$ be a Petri net. We denote by $\equi{N}$ the congruence relation on paths of $G(N)$ generated by the identification of subpaths
$$
\mu_1\red{\alpha}\nu_1\red{\beta}\mu_2 \et \mu_1\red{\beta}\nu_2\red{\alpha}\mu_2,
$$

\noindent such that the following equalities hold for a given marking $\rho$ in $\Mr(N)$ and for every place $x$ in $X$:
$$
\begin{array}{l c l c l c l}
\mu_1(x) &=&\rho(x)+w(x,\alpha)+w(x,\beta), 
&\qquad &\nu_1(x) &=& \rho(x)+w'(\alpha,x)+w(x,\beta), \\
\nu_2(x) &=&\rho(x)+w(x,\alpha)+w'(\beta,x), 
&\qquad &\mu_2(x) &=& \rho(x)+w'(\alpha,x)+w'(\beta,x).
\end{array}
$$
\end{ntt}

\noindent One can check that, in the reduction graph of the Petri net of example \ref{ex:reseau_de_petri}, the relation $\equi{N}$ identifies any two paths with same source and same target one can form in the diagram of example \ref{ex:graphe_de_reduction}. In each one of the next three sections, we introduce a translation for Petri nets and study how it behaves with respect to this congruence relation.

\section{Petri nets and commutative word rewriting systems}\label{sec:mots_commutatifs}

\noindent In this section, an equivalence between Petri nets and \emph{commutative word rewriting systems} is proved. The underlying idea of the translation is already present in [Caprotti Ferscha Hong 1995] and [Chandler Heyworth 2001] and comes from the following remarks :
\begin{enumerate}
\item[-] The markings of a Petri net have a commutative monoid structure: the sum is given by addition of the tokens in each place and the empty marking is a neutral element for this operation.
\item[-] If $\alpha$ is a transition, then $\red{\alpha}$ is compatible with the commutative monoid structure on markings: if~$\mu\red{\alpha}\mu'$, then $\mu+\nu\red{\alpha}\mu'+\nu$ holds for every marking $\nu$.
\end{enumerate}

\begin{defn}\label{def:reecriture_commutative}
Let $X$ be a set. The \emph{free commutative monoid generated by $X$} is the set $[X]$ of all finite formal sums of elements of $X$:
$$
a=\sum_{x\in X} a_x.x,
$$

\noindent where the $a_x$ are natural numbers that entirely define $a$. The set $[X]$ is a commutative monoid for the following operation, which admits the empty sum as a neutral element:
$$
\sum_{x\in X} a_x.x \:+\: \sum_{x\in X} b_x.x \:=\: \sum_{x\in X} (a_x+b_x).x.
$$

\noindent A \emph{(finite) commutative word rewriting system} is a pair $(X,R)$, where $X$ is a (finite) set, called the \emph{alphabet}, and $R$ is a (finite) family of pairs of elements of $[X]$, called the \emph{rules}. If $\alpha=(s(\alpha),t(\alpha))$ is in $R$, the \emph{reduction relation} $\red{\alpha}$ it generates is defined by $a\red{\alpha}b$ if there exists some formal sum $c$ such that $a=c+s(\alpha)$ and $b=c+t(\alpha)$. To any commutative word rewriting system $(X,R)$, one associates a \emph{reduction graph} $G(X,R)$, defined by:
\begin{enumerate}
\item[0.] The objects of $G(X,R)$ are the elements of $[X]$.
\item[1.] The arrows of $G(X,R)$ are the pairs $(c,\alpha)$ made of an element $c$ of $[X]$ and a rule $\alpha=(s(\alpha),t(\alpha))$ in $R$. Such an arrow has source $c+s(\alpha)$ and target $c+t(\alpha)$; it can be written $c+\alpha$.
\end{enumerate}
\end{defn}

\begin{rem}
The arrows of $G(X,R)$ are contextual applications of rules: indeed, there is an arrow $(c,\alpha)$ in $G(X,R)$ from $a$ to $b$ if and only if $a\red{\alpha}b$. Furthermore, in this case, $c$ is the context of the application of $\alpha$ at $a$: this is the part that remains unchanged after action of the rule.
\end{rem}

\begin{rem}
In [Guiraud 2004(T)], commutative word rewriting systems are seen as presentations by generators and relations of commutative monoids: indeed, such an object defines a commutative monoid which elements are the connected components of its reduction graph. Conversely, every commutative monoid admits a commutative word rewriting system as a presentation: the generators are the elements of the monoid and the relations are given by the "multiplication" table of the sum.
\end{rem}

\noindent Following the same idea as in section \ref{sec:reseaux_de_petri}, let us define a congruence relation between paths of the reduction graph of a commutative word rewriting system:

\begin{ntt}
Let $(X,R)$ be a commutative word rewriting system. The relation $\equi{(X,R)}$ is the congruence relation on paths of $G(X,R)$ generated by the identification of squares of the following shape, with $\alpha$ and $\beta$ in $R$ and $c$ in $[X]$: 
$$
\xymatrix
{
	c + s(\alpha) + s(\beta) \ar[rr]^-{(c+s(\beta))+\alpha} \ar[dd]_-{(c+s(\alpha))+\beta}
	&&
	c + t(\alpha) + s(\beta) \ar[dd]^-{(c+t(\alpha))+\beta} 
	\\ \\
	c + s(\alpha) + t(\beta) \ar[rr]_-{(c+t(\beta))+\alpha}
	&&
	c + t(\alpha) + t(\beta).
}
$$
\end{ntt}

\noindent Translations between Petri nets and finite commutative word rewriting systems are defined as follows:

\begin{defn}
Let $N=(X,T,w,w')$ be a Petri net. Its \emph{associated finite commutative word rewriting system} is denoted by $\Phi(N)$ and defined by:
\begin{enumerate}
\item[-] The alphabet of $\Phi(N)$ is the set $X$ of places of $N$.
\item[-] The rules of $\Phi(N)$ are the transitions of $N$, seen as pairs $\alpha=(s(\alpha),t(\alpha))$ with:
$$
s(\alpha) \:=\: \sum_{x\in X} w(x,\alpha).x \et t(\alpha) \:=\: \sum_{x\in X} w'(\alpha,x).x.
$$
\end{enumerate}

\noindent Conversely, let $(X,R)$ be a finite commutative word rewriting system. Its \emph{associated Petri net} is denoted by~$\Psi(X,R)$ and defined by:
\begin{enumerate}
\item[-] The places of $\Psi(X,R)$ are the elements of $X$.
\item[-] There is one transition in $\Psi(X,R)$ for each rule in $R$.
\item[-] The weights $w$ and $w'$ are given, on a place $x$ and a transition $\alpha=(s(\alpha),t(\alpha))$, by:
$$
w(x,\alpha) \:=\: s(\alpha)_x \et w'(\alpha,x) \:=\: t(\alpha)_x.
$$
\end{enumerate}
\end{defn}

\begin{ex}\label{ex:mots_commutatifs}
Let us consider the Petri net from example \ref{ex:reseau_de_petri}. The corresponding commutative word rewriting system is the pair $(X,R)$, where $X=\ens{x,y,z}$ and $R$ consists of the two following rewriting rules $\alpha:x\fl y+z$ and $\beta\::\:2y\fl z$. The marking from example \ref{ex:marquage} corresponds to the formal sum $2x+2y$. The reduction graph from example \ref{ex:graphe_de_reduction} becomes:
$$
\xymatrix{ 2x+2y \ar[dd]_-{x+2y+\alpha} \ar[rr]^-{2x+\beta} && 2x+z \ar[dd]^-{x+z+\alpha} \\ \\ x+3y+z \ar[dd]_-{3y+z+\alpha} \ar[rr]^-{x+y+z+\beta} && x+y+2z \ar[dd]^-{y+2z+\alpha} \\ \\ 4y+2z \ar[rr]_-{2y+2z+\beta} && 2y+3z \ar[rr]_-{3z+\beta} && 4z.}
$$

\noindent One can check that, in this diagram, any two paths with same source and same target are identified by the congruence $\equi{(X,R)}$: the translation from Petri nets to commutative word rewriting systems seems to preserve the congruence relation we have defined on Petri nets reduction paths.
\end{ex}

\noindent The following result proves that, in essence, Petri nets and finite commutative word rewriting systems are the same objects and generate the same reduction graphs:

\begin{thm}\label{th:comm}
For every Petri net $N$, the equality $\Psi\circ\Phi(N)=N$ holds and the reduction graphs~$G(N)$ and $G(\Phi(N))$ are isomorphic. Furthermore, this isomorphism identifies the congruences $\equi{N}$ and $\equi{\Phi(N)}$. Conversely, for every finite commutative word rewriting system $(X,R)$, the equality $\Phi\circ\Psi(X,R)=(X,R)$ holds and the reduction graphs $G(X,R)$ and $G(\Psi(X,R))$ are isomorphic. Furthermore, this isomorphism identifies the congruences $\equi{(X,R)}$ and $\equi{\Psi(X,R)}$.
\end{thm}

\begin{dem}
Let us fix $N=(X,T,w,w')$ and prove the equality $\Psi\circ\Phi(N)=N$. The places of $\Psi\circ\Phi(N)$ are the elements of the alphabet of $\Phi(N)$: these are the places of $N$. The transitions of $\Psi\circ\Phi(N)$ are the rules of $\Phi(N)$: these are the transitions of $N$. Let us fix a place $x$ in $X$ and a transition $\alpha$ in $T$. Let us denote by $\ol{w}$ and $\ol{w}'$ the weighting functions of $\Psi\circ\Phi(N)$ and compare them with $w$ and $w'$. By definition of $\Psi\circ\Phi(N)$:
$$
\ol{w}(x,\alpha) \:=\: s(\alpha)_x \et \ol{w}'(\alpha,x) \:=\: t(\alpha)_x.
$$

\noindent And by definition of $\Phi(N)$:
$$
s(\alpha) \:=\: \sum_{y\in X} w(y,\alpha).y \et t(\alpha) \:=\: \sum_{y\in X} w'(\alpha,y).y.
$$

\noindent Invoking the fact that $[X]$ is free, one gets:
$$
s(\alpha)_x \:=\: w(x,\alpha) \et t(\alpha)_x \:=\: w'(x,\alpha).
$$

\noindent Hence $w=\ol{w}$ and $w'=\ol{w}'$. Now, let us prove that $G(N)$ and $G(\Phi(N))$ are isomorphic graphs. We define a graph morphism $\phi$ from the former to the latter. Let $\mu$ be a marking of $N$ and let us define an element $\phi(\mu)$ in $[X]$ this way:
$$
\phi(\mu) \:=\: \sum_{x\in X} \mu(x).x.
$$

\noindent Now, let us consider an arrow $f:\mu\fl\nu$ in $G(N)$. By definition of $G(N)$, this arrow corresponds to a transition $\alpha$ such that $\mu\red{\alpha}\nu$. By definition of the relation $\red{\alpha}$ on markings, this means that:
$$
\mu(x)\geq w(x,\alpha) \et \nu(x)=\mu(x)-w(x,\alpha)+w'(\alpha,x).
$$

\noindent Let us prove that $\phi(\mu)\red{\alpha}\phi(\nu)$ is a reduction generated by $(X,R)$. By definition of $\phi$ on markings:
$$
\Phi(\mu)=\sum_{x\in X} \mu(x).x \et \Phi(\nu)=\sum_{x\in X} \nu(x).x.
$$

\noindent Hence, proving $\phi(\mu)\red{\alpha}\phi(\nu)$ is equivalent to prove that there exists a $c$ in $[X]$ such that:
$$
\sum_{x\in X}\mu(x).x=c+s(\alpha) \et \sum_{x\in X}\nu(x).x=c+t(\alpha),
$$

\noindent Since $\mu(x)\geq w(x,\alpha)$ for every place $x$, the following $c$ is well-defined in $[X]$:
$$
c=\sum_{x\in X}(\mu(x)-w(x,\alpha)).x.
$$

\noindent Then:
$$
c+s(\alpha) \:=\: \sum_{x\in X}(\mu(x)-w(x,\alpha)).x \:+\: \sum_{x\in X} w(x,\alpha).x \:=\: \sum_{x\in X} \mu(x).x.
$$

\noindent Furthermore, using the fact that $\nu(x)=\mu(x)-w(x,\alpha)+w'(\alpha,x)$ holds for every $x$, one gets:
$$
c+t(\alpha) \:=\: \sum_{x\in X}(\mu(x)-w(x,\alpha)).x \:+\: \sum_{x\in X}w'(\alpha,x).x \:=\: \sum_{x\in X}\nu(x).x.
$$

\noindent Hence $\phi(\mu)\red{\alpha}\phi(\nu)$ holds in $[X]$. By definition of $G(\Phi(N))$, this reduction corresponds to an arrow of the form $c+\alpha$, with $c$ in $[X]$, going from $\phi(\mu)$ to $\phi(\nu)$ in $G(\Phi(N))$. Let us define $\phi(f)$ to be this arrow. 

Let us define a graph morphism $\psi$ from~$G(\Phi(N))$ to~$G(N)$ and prove that it is inverse of $\phi$. Let $a$ be an element of $[X]$. Then $\psi(a)$ is defined as the marking $\psi(a)(x)=a_x$ for every place $x$. Now, let us consider an arrow $c+\alpha$ in $G(\Phi(N))$, which starts at $a=c+s(\alpha)$ and ends at~$b=c+t(\alpha)$. Then, for every place $x$:
$$
\psi(a)(x) \:=\: a_x \:=\: c_x+s(\alpha)_x \:=\: c_x+w(x,\alpha).
$$

\noindent Thus $\psi(a)(x)\geq w(x,\alpha)$. Furthermore:
$$
\psi(b)(x) \:=\: b_x \:=\: c_x+t(\alpha)_x \:=\: \psi(a)(x) - w(x,\alpha)+w'(\alpha,x).
$$

\noindent Hence $\psi(a)\red{\alpha}\psi(b)$ holds in $\Mr(N)$. This reduction corresponds to an arrow in $G(N)$, which we take as $\psi(c+\alpha)$. Checking that $\psi$ is a left and right inverse for $\phi$ is straightforward. 

In order to prove that $\phi(\equi{N})$ is $\equi{\Phi(N)}$, we prove that $\phi(\equi{N})$ is included into $\equi{\Phi(N)}$ and that $\psi(\equi{\Phi(N)})$ is included into $\equi{N}$. Furthermore, since $\phi$ and $\psi$ are graph morphisms, it is sufficient to prove these inclusions on paths of minimal lenghts, such as given in the definitions of both congruences.

Hence, let us consider two paths $\mu_1\red{\alpha}\nu_1\red{\beta}\mu_2$ and $\mu_1\red{\beta}\nu_2\red{\alpha}\mu_2$ in $G(N)$ such that there exists a marking $\rho$ of $N$ that satisfies the following four equalities for every place $x$:
$$
\begin{array}{l l l c l l l}
\mu_1(x) &=& \rho(x)+w(x,\alpha)+w(x,\beta), 
&\qquad &\nu_1(x) &=& \rho(x)+w'(\alpha,x)+w(x,\beta), \\
\nu_2(x) &=& \rho(x)+w(x,\alpha)+w'(\beta,x), 
&\qquad &\mu_2(x) &=& \rho(x)+w'(\alpha,x)+w'(\beta,x).
\end{array}
$$

\noindent Let us denote by $c$ the element $\phi(\rho)$ of $[X]$. Then $\phi$ sends both paths onto the following ones, which are identified by $\equi{\Phi(N)}$:
$$
\xymatrix
{
	c + s(\alpha) + s(\beta) \ar[rr]^-{(c+s(\beta))+\alpha} 
	&& 
	c + t(\alpha) + s(\beta) \ar[rr]^-{(c+t(\alpha))+\beta}
	&&
	c + t(\alpha) + t(\beta)
}
$$

\noindent and:
$$
\xymatrix
{
	c + s(\alpha) + s(\beta) \ar[rr]^-{(c+s(\alpha))+\beta} 
	&& 
	c + s(\alpha) + t(\beta) \ar[rr]^-{(c+t(\beta))+\alpha}
	&&
	c + t(\alpha) + t(\beta).
}
$$

\noindent Then, let us consider two paths in $G(\Phi(N))$ written as above, for a given $c$ in $[X]$. Let us denote by $\rho$ the marking $\psi(c)$. Then, if the four markings $\mu_1$, $\mu_2$, $\nu_1$ and $\nu_2$ are defined as above, the graph morphism $\psi$ sends both paths of $G(\Phi(N))$ onto $\mu_1\red{\alpha}\nu_1\red{\beta}\mu_2$ and $\mu_1\red{\beta}\nu_2\red{\alpha}\mu_2$: these two paths are identified by~$\equi{N}$. 

Conversely, let us consider a finite commutative word rewriting system $(X,R)$ and prove that the equality $\Phi\circ\Psi(X,R)=(X,R)$ holds. By definition of the rewriting system $\Phi\circ\Psi(X,R)$, its alphabet is the set of places of $\Psi(X,R)$: this is the alphabet of $(X,R)$. The rules in $\Phi\circ\Psi(X,R)$ are the pairs $(\ol{s}(\alpha),\ol{t}(\alpha))$ for each transition $\alpha$ in $\Psi(X,R)$, where:
$$
\ol{s}(\alpha) \:=\: \sum_{x\in X} w(x,\alpha).x \et \ol{t}(\alpha) \:=\: \sum_{x\in X}w'(\alpha,x).x.
$$

\noindent Furthermore, each transition $\alpha$ in $\Psi(X,R)$ comes from a rule $(s(\alpha),t(\alpha))$ in $R$ and:
$$
w(x,\alpha) \:=\: s(\alpha)_x \et w'(\alpha,x) \:=\: t(\alpha)_x.
$$

\noindent Thus, $\ol{s}(\alpha)=s(\alpha)$ and $\ol{t}(\alpha)=t(\alpha)$, so that the set of rules of $\Phi\circ\Psi(X,R)$ is $R$. Hence, the two commutative word rewriting systems $(X,R)$ and $\Phi\circ\Psi(X,R)$ are the same. 

Let us prove that the two graphs $G(X,R)$ and $G(\Psi(X,R))$ are isomorphic. Since $\Psi(X,R)$ is a Petri net, we already know that $G(\Psi(X,R))$ is isomorphic to $G(\Phi\circ\Psi(X,R))$: this graph is $G(X,R)$ since the equality $\Phi\circ\Psi(X,R)=(X,R)$ holds. Furthermore, this graph isomorphism is defined the same way as $\phi$ and $\psi$ in the first part of the proof. Hence $\phi(\equi{\Psi(X,R)})$ is equal to $\equi{(X,R)}$. If one applies $\psi$, one gets the equality of both congruences $\equi{\Psi(X,R)}$ and $\psi(\equi{(X,R)})$.

\findem\end{dem}

\begin{rem}
This equivalence between Petri nets and finite commutative word rewriting systems highlights the underlying algebraic structure of the formers: one immediate usage is that every arrow in the reduction graph has an explicit name, such as $x+2y+\alpha$, giving the context of application of the rule~$\alpha$. 
\end{rem}

\begin{rem}
Another more concrete concrete usage of the translation was developped in the aforementioned [Caprotti Ferscha Hong 1995] and [Chandler Heyworth 2001]: there, it was decribed how Gr\"obner bases can be used to partially solve the reachability problem for Petri nets, when they are seen as commutative word rewriting systems. 
\end{rem}

\begin{rem}
If $N$ is a Petri net, the definition of $\equi{N}$ is technical but intuitively simple. The unveiling of the intrinsic algebraic structure of Petri nets makes this definition much simpler. Indeed, let us consider a commutative word rewriting system $(X,R)$ and denote by $\circ$ the composition of paths in the graph $G(X,R)$. Note that this amounts at considering the category $\mon{G(X,R)}$ freely generated by $G(X,R)$, as it is defined in section \ref{sec:2d}. Then, the relation $\equi{(X,R)}$ can be defined as the congruence on $\mon{G(X,R)}$ generated by the following identifications, for any $c$ in $[X]$:
$$
(c+t(\alpha)+\beta)\circ(c+s(\beta)+\alpha) \:\equiv\: (c+t(\beta)+\alpha)\circ(c+s(\alpha)+\beta).
$$

\noindent Let us also note that such equations allow the sum of $[X]$ to be naturally extended to reductions: $\alpha+\beta$ will be any side of the given equation for $c=0$. This is also the idea developped with polygraphs in sections~\ref{sec:2d} and~\ref{sec:3d}.
\end{rem}

\noindent From now on, theorem \ref{th:comm} grants us the right to consider that a Petri net \emph{is} a finite commutative word rewriting system. In fact, the results to be proved are not limited to the finite case. Hence, thereafter, the name Petri net stands for a commutative word rewriting system. Let us use this new equivalent definition to give a different graphical account of Petri nets.

\section{Petri nets as $\mathbf{2}$-dimensional objects}\label{sec:2d}

\noindent The goal of this section is to prove that Petri nets have strong links with a certain class of \emph{two-dimensional polygraphs}. The first result presented here, theorem~\ref{th:2d}, is essentially due to Albert Burroni, who gived the idea of the translation. The behaviour of this translation with respect to the congruence on Petri nets reduction paths is new and described in proposition~\ref{prop:2d-equiv}. A discussion follows on many issues to be studied in future work.

In order to translate Petri nets into polygraphs, we start by the interpretation of the markings of a Petri net (the formal sums of its places) into $1$-dimensional objects. Let us recall the some classical notions about graphs, free categories and monoids.

\begin{ntt}
If $G$ is a graph, its set of objects is denoted by $G_0$ and its set of arrows going from an object~$x$ to another object $y$ is denoted by $G(x,y)$; for such an arrow $f$, $s_0(f)$ is the source $x$ of $f$ and~$t_0(f)$ its target $y$. The set of all arrows of $G$ is denoted by $G_1$ and $G$ itself is often abusively denoted by~$(G_0,G_1)$ only, assuming that the source and target mappings are given with $G_1$. 
\end{ntt}

\begin{defn}
Let $G=(G_0,G_1)$ be a graph. The \emph{free category generated by $G$}, denoted by $\mon{G}$, is the following (small) category:
\begin{enumerate}
\item[0.] The objects of $\mon{G}$ are the objects of $G$.
\item[1.] The arrows of $\mon{G}$, from $x$ to $y$, are all the finite paths in $G$ going from $x$ to $y$. Their composition, denoted by $\circ$, is the concatenation of paths. The empty paths are local identities for this operation.
\end{enumerate}

\noindent Such a category is often denoted by $\mon{G}=(\mon{G}_0,\mon{G}_1)$ or just by $(G_0,\mon{G}_1)$, assuming that the source and target mappings are given with the data in $\mon{G}_1$, together with the identities and composition operations. 
\end{defn}

\begin{ex}
Let $G=(\ast,X)$ be a graph with only one object ($\ast$ denotes any single-element set); the set of arrows can be any set $X$, with source and target being the only possible map from $X$ to $\ast$. Then the free category $\mon{G}$ is the free monoid $\mon{X}$ generated by $X$: more precisely, the set $\mon{G}(\ast,\ast)$, containing all the arrows of $\mon{G}$, equipped with the composition and the identity of $\ast$, is isomorphic to the free monoid~$\mon{X}$. A proof can be found in [MacLane 1998], for example.
\end{ex}

\begin{defn}
Let $\Cr$ be a category. Two arrows in $\Cr$ are \emph{parallel} when they have same source and same target. A \emph{relation} in $\Cr$ is a pair of parallel arrows of $\Cr$. If $R$ is a family of relations in $\Cr$, the \emph{quotient of~$\Cr$ by $R$} is the category denoted by $\Cr/R$ built this way:
\begin{enumerate}
\item[0.] The objects of $\Cr/R$ are the objects of $\Cr$.
\item[1.] The arrows from $x$ to $y$ in $\Cr/R$ are the elements of $\Cr(x,y)$, \emph{modulo} the reflexive-symmetric-transitive closure $\equi{R}$ of the relation $\red{R}$ defined by: $f\red{R}g$ if there exist a relation $(u,v)$ in $R$ and two arrows $h$ and $k$ in $\Cr$ such that $f=k\circ u\circ h$ and $g=k\circ v\circ h$. The identities of $\Cr/R$ are the equivalence classes of the identity of $\Cr$. The composition of $\Cr/R$ is induced by the one in $\Cr$.
\end{enumerate}
\end{defn}

\begin{rem}
The defined object $\Cr/R$ is only a graph. One must check, through easy computations, that the composition of $\Cr$ is compatible with $\equi{R}$: the result of the composition is independent of any choice of representatives. Furthermore, it must be checked that induced composition satisfies the axioms of associativity and left and right units of the category structure.
\end{rem}

\begin{ex}
Let $G=(\ast,X)$ be a graph with one object. On $\mon{G}$, one defines $R$ to be the family of all relations $(x\circ y,y\circ x)$, for $x$ and $y$ in $X$. Then $\mon{G}/R$ is the free commutative monoid $[X]$ generated by the set $X$.
\end{ex}

\noindent Hence, we have a graphical description of $[X]$. However, the main idea behind higher-dimensional rewriting is to replace any equation between $n$-dimensional objects by a $(n+1)$-dimensional object: equalities are replaced by their proofs - this point of view was developped in both [Burroni 1993] and [Baez Dolan 1998]. Following this leading idea, equalities of the form $x\circ y=y\circ x$ are replaced by $2$-dimensional cells, pasted between parallel paths in the graph $(\ast,X)$, such as the following one:
$$
\xymatrix{& \ast \ar[dr]^y \ar@{=>}-<0cm,0.5cm>;[dd]+<0cm,0.5cm>^{\tau_{x,y}} \\ \ast \ar[ur]^x \ar[dr]_y && \ast \\ & \ast \ar[ur]_x}
$$

\begin{rem}
In order to achieve commutativity, one may ask that $\tau_{x,y}$ is an isomorphism, with $\tau_{y,x}$ as inverse: in this case, one gets a \emph{categorified} version of the free commutative monoid. Another point of view would be to replace the equalities $\tau_{y,x}\circ\tau_{x,y}=\id_{x\tens y}$ and $\tau_{x,y}\circ\tau_{y,x}=\id_{y\tens x}$ by their proofs: these would be $3$-dimensional cells. This issue is discussed at the end of this section. 
\end{rem}

\noindent So far, we have described an object with one $0$-cell, as many $1$-cells as there are in our set $X$, together with one $2$-cell $\tau_{x,y}$ for each pair $(x,y)$ of distinct elements in $X$. Now, let us consider the rule $\alpha:x\fl y+z$ from example \ref{ex:mots_commutatifs}. Such a rule is also translated as a $2$-dimensional cell:
$$
\xymatrix{\ast \ar[rr]^x_{}="1" \ar@{=>}"1"-<0cm,0.15cm>;"2,2"+<0cm,0.3cm>^{\alpha} \ar[dr]_y && \ast \\ & \ast \ar[ur]_z }
$$

\noindent A choice has been made in order to represent the rule $\alpha$. Indeed, it could have been seen as transforming~$x$ into $z+y$, which is equal to $y+z$ in the commutative monoid $[X]$. This is the arbitrary part of the presented $2$-polygraphic interpretation of Petri nets: it assumes that, for every element $a$ in $[X]$, a representative has been chosen in $\mon{X}$. 

Since we must use the axiom of choice, let us apply the equivalent Zermelo theorem and assume, until the end of this section, that, for every Petri net $(X,R)$, the set $X$ comes equipped with a total order. Then, every element $a$ of $[X]$ has a unique decomposition $a=n_1.x_1+\dots+n_k.x_k$, where the $n_i$ are non-zero natural numbers and the $x_i$ are elements of $X$ such that $x_1<\dots<x_k$. 

\begin{ntt}
Let $X$ be a set and $a$ an element of $[X]$. Let us denote by $n_1.x_1+\dots+n_k.x_k$ the unique decomposition of $a$. Then $\ol{a}$ denotes the representative $x_1^{n_1}\dots x_k^{n_k}$ of $a$ in $\mon{X}$, where $x^n$ is the product in $\mon{X}$ of $n$ copies of $x$.
\end{ntt}

\noindent Until now, we have constructed a composite object $\Sigma=(\Sigma_0,\Sigma_1,\Sigma_2)$, made of sets $\Sigma_i$ of $i$-dimensional cells. On top of these three sets, $\Sigma$ also contains boundaries informations: for example, the cell~$\tau_{x,y}$ has source~$x\circ y$ and target $y\circ x$, while $\alpha$ has source $x$ and target $y\circ z$. 

Such an object is called a \emph{polygraph}: it is the central structure studied in \emph{higher-dimensional rewriting}. Here, the object $\Sigma$ is a $2$-dimensional polygraph or $2$-polygraph for short. Its definition is recalled from [Burroni 1993].

\begin{defn}\label{def:2-polygraphe}
A \emph{$\mathit{2}$-polygraph} $\Sigma$ is given by:
\begin{enumerate}
\item[0.] A set $\Sigma_0$ of \emph{$\mathit{0}$-cells}.
\item[1.] A set $\Sigma_1$ of \emph{$\mathit{1}$-cells}, together with two maps $s_0,t_0:\Sigma_1\fl\Sigma_0$, called \emph{$\mathit{0}$-source} and \emph{$\mathit{0}$-target}. The arrows of the free category $(\Sigma_0,\mon{\Sigma}_1)$ are called \emph{$\mathit{1}$-arrows}. The composition of $f$ followed by $g$ is denoted by $f\star_0g$ or $g\circ_0f$ in the general case and $f\tens g$ when $\Sigma_0$ has only one element.
\item[2.] A set $\Sigma_2$ of \emph{$\mathit{2}$-cells}, together with two maps $s_1,t_1:\Sigma_2\fl\mon{\Sigma}_1$, called \emph{$\mathit{1}$-source} and \emph{$\mathit{1}$-target}, and such that $s_0\circ s_1=s_0\circ t_1$ and $t_0\circ s_1 = t_0\circ t_1$. The first equality gives a map $s_0:\Sigma_2\fl\Sigma_0$ and the second one yields $t_0:\Sigma_2\fl\Sigma_0$.
\end{enumerate}
\end{defn}

\begin{defn}
Let $(X,R)$ be a Petri net, such that $X$ is equipped with a total order. The \emph{$\mathit{2}$-polygraph associated with $(X,R)$} is $\Sigma^2(X,R)$ defined this way:
\begin{enumerate}
\item[0.] There is one $0$-cell in $\Sigma^2(X,R)$, denoted by $\ast$.
\item[1.] The $1$-cells of $\Sigma^2(X,R)$ are the elements of $X$, with the only possible $0$-source and $0$-target maps.
\item[2.] The $2$-cells of $\Sigma^2(X,R)$ consist of all the $\tau_{x,y}$, for $x\neq y$ in $X$, together with one $2$-cell $\alpha$ for each rule in $R$. The $1$-source and $1$-target maps are given by:
$$
s_1(\tau_{x,y})\:=\: x\tens y, \quad t_1(\tau_{x,y})\:=\: y\tens x, \quad s_1(\alpha)\:=\: \ol{s(\alpha)}, \quad t_1(\alpha)\:=\: \ol{t(\alpha)}.
$$
\end{enumerate}
\end{defn}

\noindent In order to compare a Petri net to its associated $2$-polygraph, we define a notion of reduction graph for these objects. The idea is to see every $2$-cell of a $2$-polygraph as a rewriting rule on $1$-arrow, that can be applied in any context: a $2$-cell $\alpha$ can be applied on any $1$-arrow of the shape $u\tens s_1(\phi)\tens v$, in order to produce the $1$-arrow $u\tens t_1(\phi)\tens v$. Let us formalize this idea.

\begin{defn}
Let $\Sigma=(\Sigma_0,\Sigma_1,\Sigma_2)$ be a $2$-polygraph. The \emph{reduction graph associated to $\Sigma$}, denoted by $G(\Sigma)$, is defined this way:
\begin{enumerate}
\item[0.] The objects of $G(\Sigma)$ are the $1$-arrows of $\Sigma$. 
\item[1.] The arrows from $f$ to $g$ in $G(\Sigma)$ are the triples $(h,\phi,k)$ where $h$ and $k$ are $1$-arrows in $\mon{\Sigma}_1$ and~$\phi$ is a $2$-cell in $\Sigma_2$ such that the following equalities hold:
$$
f \:=\: h\star_0 s_1(\phi)\star_0 k \et g \:=\: h\star_0 t_1(\phi)\star_0 k.
$$

\noindent A triple $(h,\phi,k)$ is denoted $h\star_0\phi\star_0 k$, and $h\star_0$ (resp. $\star_0 k$) is dropped when $h$ (resp. $k$) is an identity (an empty path).
\end{enumerate}
\end{defn}

\noindent We want to prove that the two graphs $G(X,R)$ and $G(\Sigma^2(X,R))$ have strong links. To begin with, let us note that the objects of the graph $G(\Sigma^2(X,R))$ are the elements of the free monoid $\mon{X}$, while the objects of the graph $G(X,R)$ are the ones of the free commutative monoid $[X]$. We define $\pi:\mon{X}\mfl[X]$ to be the canonical projection. 

\begin{lem}\label{lem:relevement}
Let $u$ and $v$ be two elements in $\mon{X}$ such that $\pi(u)=\pi(v)$. Then, there exists an arrow $f$ in $G(\Sigma^2(X,R))$ with source $u$ and target $v$, such that $f$ has a decomposition of the form:
$$
f \:=\: (u_n\tens\tau_{x_n,y_n}\tens v_n)\circ\dots\circ(u_1\tens\tau_{x_1,y_1}\tens v_1).
$$
\end{lem}

\begin{dem}
Since $\mon{X}$ is freely generated by $X$, the elements $u$ and $v$ uniquely decompose as:
$$
u \:=\: z_1\tens\dots\tens z_p \et v \:=\: z'_1\tens\dots\tens z'_{p'},
$$

\noindent with the $z_i$ and $z'_i$ in $X$. Since $\pi(u)=\pi(v)$, the following equality holds in $[X]$:
$$
z_1+\dots+z_p \:=\: z'_1+\dots+z'_{p'}.
$$

\noindent Hence, since $[X]$ is freely generated by $X$, we get that $p=p'$ and that there exists a permutation $\sigma$ in $\Sk_p$ such that, for every $i\in\dens{1}{p}$, $z'_{\sigma(i)}=z_i$. Let us consider a decomposition of the permutation $\sigma$ in~$n$ transpositions:
$$
\sigma \:=\: \tau_{i_n}\circ\dots\circ\tau_{i_1},
$$

\noindent where each $i_j$ is in $\dens{1}{p-1}$ and $\tau_{i_j}$ is the transposition that exchanges $i_j$ and $i_{j+1}$. Let us fix the following notations:
$$
u_1 \:=\: z_1\tens\dots\tens z_{i_1-1}, \quad x_1 \:=\: z_{i_1}, \quad y_1 \:=\: z_{i_1+1}, \quad v_1 \:=\: z_{i_1+2}\tens\dots\tens z_p.
$$

\noindent Then, the arrow $f_1=u_1\tens\tau_{x_1,y_1}\tens v_1$ of $G(\Sigma^2(X,R))$ has source $u$ and target:
$$
z_1\tens\dots\tens z_{i_1-1}\tens z_{i_1+1}\tens z_{i_1}\tens z_{i_1+2}\tens\dots\tens z_p.
$$

\noindent But this element of $\mon{X}$ can also be written as $z_{\tau_{i_1}(1)}\tens\dots\tens z_{\tau_{i_1}(p)}$. Hence, if we repeat this construction for each $\tau_{i_j}$, we prove, by induction on the length of the decomposition of $\sigma$, that the target of the last arrow $f_n=u_n\tens\tau_{x_n,y_n}\tens v_n$, associated with $\tau_{i_n}$, is:
$$
v \:=\: z_{\sigma(1)}\tens\dots\tens z_{\sigma(p)}.
$$

\noindent In conclusion, $f=f_n\circ\dots\circ f_1$ satisfies the required hypotheses.

\findem\end{dem}

\noindent Now, the main result of this section can be proved. As mentioned earlier, this result formalizes a construction due to Albert Burroni:

\begin{thm}\label{th:2d}
Let $(X,R)$ be a Petri net. The following equalities extend the canonical map $\pi$ into a surjective functor from the free category~$\mon{G(\Sigma^2(X,R))}$ to the free category $\mon{G(X,R)}:$
$$
\pi(u\tens\tau_{x,y}\tens v) \:=\: \id_{\pi(u)+x+y+\pi(v)} \et \pi(u\tens\alpha\tens v) \:=\: \pi(u)+\pi(v)+\alpha.
$$
\end{thm}

\begin{dem}
The equalities extend $\pi$ so that it is now defined on every object and arrow of the reduction graph~$G(\Sigma^2(X,R))$ and takes its values into the free category $\mon{G(X,R)}$. Hence, a classical categorical argument tells us that $\pi$ uniquely extends into a functor, still denoted by $\pi$, from the free category~$\mon{G(\Sigma^2(X,R))}$ to the free category $\mon{G(X,R)}$. Now, let us prove that $\pi$ is surjective, which means that both its restrictions on objects and on arrows are surjective. On objects, $\pi$ is the canonical morphism from the free monoid $\mon{X}$ to the free commutative monoid $[X]$, which is surjective. 

Let us consider two objects $a$ and $b$ in $\mon{G(X,R)}$: they are elements of the free commutative monoid~$[X]$. Let $f$ be an arrow in $\mon{G(X,R)}$ from $a$ to $b$. By definition of $G(X,R)$ and of the free category it generates, this means that $f$ uniquely decomposes as:
$$
f=(c_k+\alpha_k)\circ\dots\circ(c_1+\alpha_1),
$$

\noindent with the $c_i$ in $[X]$ and the $\alpha_i$ in $R$, such that the following relations hold in $[X]$:
$$
c_1+s(\alpha_1)\:=\: a, \quad c_i+t(\alpha_i)\:=\: c_{i+1}+s(\alpha_{i+1}), \quad c_k+t(\alpha_k)\:=\:b.
$$

\noindent Let us denote by $f_i$ the arrow $\ol{c_i}\tens\alpha_i$ in $G(\Sigma^2(X,R))$: it has source $\ol{c_i}\tens\ol{s(\alpha_i)}$ and target $\ol{c_i}\tens\ol{t(\alpha_i)}$. Hence, the equalities $\pi(s(f_1))=a$ and $\pi(t(f_n))=b$ hold. There remains to link all the $f_i$ in order to conclude. Indeed, the relation $t(f_i)=s(f_{i+1})$ does not necessarily hold for every $i$, so that $f_i$ and $f_{i+1}$ are not composable in general. 

However, the relation $\pi(t(f_i))=\pi(s(f_{i+1}))$ holds, by assumption, for every $i$. By application of lemma~\ref{lem:relevement}, we know that there exist arrows $g_1$, $\dots$, $g_{k-1}$ in $\mon{G(\Sigma^2(X,R))}$ such that each one is a composition of arrows of the form $(u\tens\tau_{x,y}\tens v)$ and such that the following diagram is an arrow of~$\mon{G(\Sigma^2(X,R))}$:
$$
\xymatrix{ \ol{a} \ar[r]^-{f_1} & \ol{c_1}\tens\ol{t(f_1)} \ar[r]^-{g_1} & \ol{c_2}\tens\ol{s(f_2)} \ar[r]^-{f_2} & \dots \ar[r]^-{g_{k-1}} & \ol{c_k}\tens\ol{s(f_k)} \ar[r]^-{f_k} & \ol{b}.}
$$

\noindent Finally, from the definition of the functor $\pi$, we conclude that:
$$
\pi(g_i)=\id_{c_i+t(f_i)} \:=\: \id_{c_{i+1}+s(f_{i+1})} \et \pi(f_i) \:=\: c_i+\alpha_i.
$$

\noindent Hence $\pi(f_k\circ g_k\circ\dots\circ g_1\circ f_1)=f$, so that $\pi$ is a surjective functor.

\findem\end{dem}

\noindent So far, we have built a new graphical object $G(\Sigma^2(X,R))$ in which every path represents a possible evolution of the Petri net $(X,R)$ and in which every possible evolution has a representative. 

But $G(\Sigma^2(X,R))$ is not the natural object one would build from the $2$-polygraph $\Sigma^2(X,R)$: indeed, such a polygraph is a presentation of a $2$-category, which is a quotient of $\mon{G(\Sigma^2(X,R))}$ by some topology-flavoured relations. Furthermore, we will see that these relations are the ones that identify the intuitively equal paths from examples \ref{ex:graphe_de_reduction} and \ref{ex:mots_commutatifs}.

Here we only define the notion of \emph{free $2$-category} generated by a $2$-polygraph with one $0$-cell, while the complete construction is in [Burroni 1993] and [Métayer 2003]. After the formal algebraic definition, we give the topological intuition that underlies it.

\begin{defn}
Let $\Sigma=(\ast,\Sigma_1,\Sigma_2)$ be a $2$-polygraph with one $0$-cell. The \emph{free $\mathit{2}$-category} generated by $\Sigma$, denoted by $\mon{\Sigma}$, is the following $2$-polygraph:
\begin{enumerate}
\item[0.] It has one $0$-cell.
\item[1.] Its $1$-cells are the $1$-arrows of $\Sigma$, which are the elements of $\mon{\Sigma}_1$.
\item[2.] Its $2$-cells, called \emph{$\mathit{2}$-arrows}, from $u$ to $v$ are the paths in the reduction graph $G(\Sigma)$, \emph{modulo} the congruence $\equi{01}$ generated by the following \emph{exchange relations} (where $g\circ f$ is written with $f$ on top of $g$ in order to match the graphical representations to be introduced):
$$
\begin{array}{c c c}
u\tens\phi\tens(v\tens s_1(\psi)\tens w) 
&& (u\tens s_1(\phi)\tens v)\tens\psi\tens w \\
\circ & \qquad\equiv\qquad & \circ \\
(u\tens t_1(\phi)\tens v)\tens\psi\tens w 
&& u\tens\phi\tens(v\tens t_1(\psi)\tens w)
\end{array}
$$

\noindent for every $2$-cells $\phi$ and $\psi$, every $1$-arrows $u$, $v$ and $w$ and where $\circ$ denotes the composition of paths in $G(\Sigma)$.
\end{enumerate}

\noindent The $2$-arrows, collectively denoted by $\mon{\Sigma}_2$, are equipped with two compositions: the first one is $\circ$, the operation yielded by the composition of paths in $G(\Sigma)$; the second one is an extension of $\tens$, allowed by the exchange relations, which is defined by functorial extension of:
$$
\begin{array}{c c c}
&& (u\tens s_1(\phi)\tens v\tens u')\tens\phi'\tens v' \\
(u\tens\phi\tens v)\tens(u'\tens\phi'\tens v') & 
\quad=\quad & \circ \\
&& u\tens\phi\tens(v\tens u'\tens t_1(\phi')\tens v')
\end{array}
$$

\end{defn}

\begin{rem}
This definition can be quite obscure and the $2$-arrows of the free $2$-category are hard to represent with the traditional cellular graphical representation. However, they become really easy to handle when using a dual representation, making the $2$-dimensional arrows appear as circuits. Let us explain how this representation is built in the case of a $2$-polygraph $\Sigma=(\ast,\Sigma_1,\Sigma_2)$ with one $0$-cell. 

Each $1$-cell $x$ is drawn as a vertical wire, labelled with $x$ (or with any symbol or color associated to the $1$-cell~$x$). A $1$-arrow is drawn as the horizontal juxtaposition of the wires representing the $1$-cells it is made of. Hence, the empty path $\id_{\ast}$ is pictured as an empty diagram and the $1$-arrow $x_1\tens\dots\tens x_n$ as:
\begin{center}
\begin{picture}(0,0)%
\includegraphics{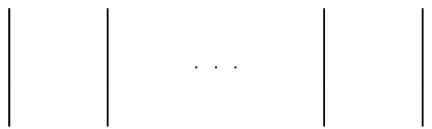}%
\end{picture}%
\setlength{\unitlength}{4144sp}%
\begingroup\makeatletter\ifx\SetFigFont\undefined%
\gdef\SetFigFont#1#2#3#4#5{%
  \reset@font\fontsize{#1}{#2pt}%
  \fontfamily{#3}\fontseries{#4}\fontshape{#5}%
  \selectfont}%
\fi\endgroup%
\begin{picture}(2350,798)(-409,107)
\put(-179,749){\makebox(0,0)[b]{\smash{{\SetFigFont{12}{14.4}{\rmdefault}{\mddefault}{\updefault}{\color[rgb]{0,0,0}$x_1$}%
}}}}
\put(1711,749){\makebox(0,0)[b]{\smash{{\SetFigFont{12}{14.4}{\rmdefault}{\mddefault}{\updefault}{\color[rgb]{0,0,0}$x_n$}%
}}}}
\put(1261,749){\makebox(0,0)[b]{\smash{{\SetFigFont{12}{14.4}{\rmdefault}{\mddefault}{\updefault}{\color[rgb]{0,0,0}$x_{n-1}$}%
}}}}
\put(271,749){\makebox(0,0)[b]{\smash{{\SetFigFont{12}{14.4}{\rmdefault}{\mddefault}{\updefault}{\color[rgb]{0,0,0}$x_2$}%
}}}}
\end{picture}%
\end{center}

\noindent A $2$-cell $\phi:u\fl v$ is pictured as a circuit component, with the wires corresponding to $u$ on top, the ones for $v$ at the bottom, such as:
\begin{center}
\begin{picture}(0,0)%
\includegraphics{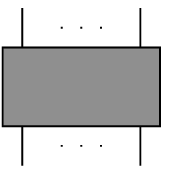}%
\end{picture}%
\setlength{\unitlength}{4144sp}%
\begingroup\makeatletter\ifx\SetFigFont\undefined%
\gdef\SetFigFont#1#2#3#4#5{%
  \reset@font\fontsize{#1}{#2pt}%
  \fontfamily{#3}\fontseries{#4}\fontshape{#5}%
  \selectfont}%
\fi\endgroup%
\begin{picture}(744,1081)(439,-335)
\put(811,-286){\makebox(0,0)[b]{\smash{{\SetFigFont{12}{14.4}{\rmdefault}{\mddefault}{\updefault}{\color[rgb]{0,0,0}$v$}%
}}}}
\put(811,164){\makebox(0,0)[b]{\smash{{\SetFigFont{12}{14.4}{\rmdefault}{\mddefault}{\updefault}{\color[rgb]{0,0,0}$\phi$}%
}}}}
\put(811,614){\makebox(0,0)[b]{\smash{{\SetFigFont{12}{14.4}{\rmdefault}{\mddefault}{\updefault}{\color[rgb]{0,0,0}$u$}%
}}}}
\end{picture}%
\end{center}

\noindent A $2$-arrow is pictured as a circuit built from the circuit components corresponding to the $2$-cells it is made of. The two compositions $\tens$ and $\circ$ are respectively represented as horizontal juxtaposition and vertical branching:
\begin{center}
\begin{picture}(0,0)%
\includegraphics{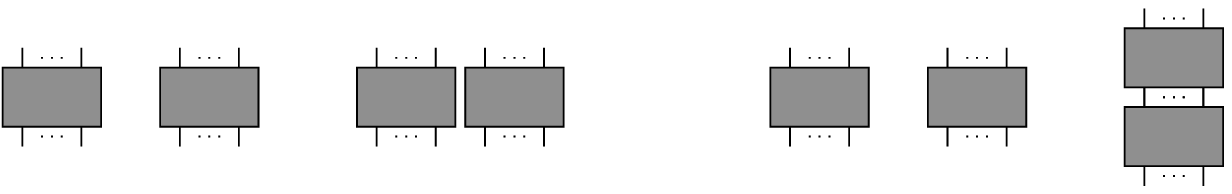}%
\end{picture}%
\setlength{\unitlength}{4144sp}%
\begingroup\makeatletter\ifx\SetFigFont\undefined%
\gdef\SetFigFont#1#2#3#4#5{%
  \reset@font\fontsize{#1}{#2pt}%
  \fontfamily{#3}\fontseries{#4}\fontshape{#5}%
  \selectfont}%
\fi\endgroup%
\begin{picture}(5604,1138)(79,-236)
\put(5446,794){\makebox(0,0)[b]{\smash{{\SetFigFont{10}{12.0}{\rmdefault}{\mddefault}{\updefault}{\color[rgb]{0,0,0}$u$}%
}}}}
\put(676,299){\makebox(0,0)[b]{\smash{{\SetFigFont{10}{12.0}{\rmdefault}{\mddefault}{\updefault}{\color[rgb]{0,0,0}$\tens$}%
}}}}
\put(316,299){\makebox(0,0)[b]{\smash{{\SetFigFont{10}{12.0}{\rmdefault}{\mddefault}{\updefault}{\color[rgb]{0,0,0}$\phi$}%
}}}}
\put(1036,299){\makebox(0,0)[b]{\smash{{\SetFigFont{10}{12.0}{\rmdefault}{\mddefault}{\updefault}{\color[rgb]{0,0,0}$\phi'$}%
}}}}
\put(316,614){\makebox(0,0)[b]{\smash{{\SetFigFont{10}{12.0}{\rmdefault}{\mddefault}{\updefault}{\color[rgb]{0,0,0}$u$}%
}}}}
\put(316,-61){\makebox(0,0)[b]{\smash{{\SetFigFont{10}{12.0}{\rmdefault}{\mddefault}{\updefault}{\color[rgb]{0,0,0}$v$}%
}}}}
\put(1036,614){\makebox(0,0)[b]{\smash{{\SetFigFont{10}{12.0}{\rmdefault}{\mddefault}{\updefault}{\color[rgb]{0,0,0}$u'$}%
}}}}
\put(1036,-61){\makebox(0,0)[b]{\smash{{\SetFigFont{10}{12.0}{\rmdefault}{\mddefault}{\updefault}{\color[rgb]{0,0,0}$v'$}%
}}}}
\put(1486,299){\makebox(0,0)[b]{\smash{{\SetFigFont{10}{12.0}{\rmdefault}{\mddefault}{\updefault}{\color[rgb]{0,0,0}$=$}%
}}}}
\put(2206,614){\makebox(0,0)[b]{\smash{{\SetFigFont{10}{12.0}{\rmdefault}{\mddefault}{\updefault}{\color[rgb]{0,0,0}$u\tens u'$}%
}}}}
\put(2206,-16){\makebox(0,0)[b]{\smash{{\SetFigFont{10}{12.0}{\rmdefault}{\mddefault}{\updefault}{\color[rgb]{0,0,0}$v\tens v'$}%
}}}}
\put(1936,299){\makebox(0,0)[b]{\smash{{\SetFigFont{10}{12.0}{\rmdefault}{\mddefault}{\updefault}{\color[rgb]{0,0,0}$\phi$}%
}}}}
\put(2431,299){\makebox(0,0)[b]{\smash{{\SetFigFont{10}{12.0}{\rmdefault}{\mddefault}{\updefault}{\color[rgb]{0,0,0}$\phi'$}%
}}}}
\put(5446,479){\makebox(0,0)[b]{\smash{{\SetFigFont{10}{12.0}{\rmdefault}{\mddefault}{\updefault}{\color[rgb]{0,0,0}$\phi$}%
}}}}
\put(5446,119){\makebox(0,0)[b]{\smash{{\SetFigFont{10}{12.0}{\rmdefault}{\mddefault}{\updefault}{\color[rgb]{0,0,0}$\psi$}%
}}}}
\put(4546,299){\makebox(0,0)[b]{\smash{{\SetFigFont{10}{12.0}{\rmdefault}{\mddefault}{\updefault}{\color[rgb]{0,0,0}$\phi$}%
}}}}
\put(3826,299){\makebox(0,0)[b]{\smash{{\SetFigFont{10}{12.0}{\rmdefault}{\mddefault}{\updefault}{\color[rgb]{0,0,0}$\psi$}%
}}}}
\put(4996,299){\makebox(0,0)[b]{\smash{{\SetFigFont{10}{12.0}{\rmdefault}{\mddefault}{\updefault}{\color[rgb]{0,0,0}$=$}%
}}}}
\put(4186,299){\makebox(0,0)[b]{\smash{{\SetFigFont{10}{12.0}{\rmdefault}{\mddefault}{\updefault}{\color[rgb]{0,0,0}$\circ$}%
}}}}
\put(4546,614){\makebox(0,0)[b]{\smash{{\SetFigFont{10}{12.0}{\rmdefault}{\mddefault}{\updefault}{\color[rgb]{0,0,0}$u$}%
}}}}
\put(4546,-16){\makebox(0,0)[b]{\smash{{\SetFigFont{10}{12.0}{\rmdefault}{\mddefault}{\updefault}{\color[rgb]{0,0,0}$v$}%
}}}}
\put(3826,614){\makebox(0,0)[b]{\smash{{\SetFigFont{10}{12.0}{\rmdefault}{\mddefault}{\updefault}{\color[rgb]{0,0,0}$v$}%
}}}}
\put(3826,-16){\makebox(0,0)[b]{\smash{{\SetFigFont{10}{12.0}{\rmdefault}{\mddefault}{\updefault}{\color[rgb]{0,0,0}$w$}%
}}}}
\put(5446,-196){\makebox(0,0)[b]{\smash{{\SetFigFont{10}{12.0}{\rmdefault}{\mddefault}{\updefault}{\color[rgb]{0,0,0}$w$}%
}}}}
\end{picture}%
\end{center}

\noindent The circuits are identified \emph{modulo} homeomorphic deformation, which exactly corresponds to the equations of the $2$-category structure. For example, the exchange relations are pictured this way:
\begin{center}
\begin{picture}(0,0)%
\includegraphics{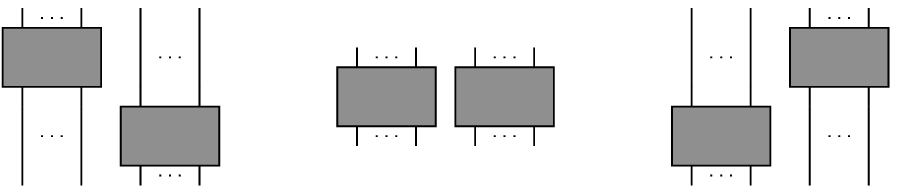}%
\end{picture}%
\setlength{\unitlength}{4144sp}%
\begingroup\makeatletter\ifx\SetFigFont\undefined%
\gdef\SetFigFont#1#2#3#4#5{%
  \reset@font\fontsize{#1}{#2pt}%
  \fontfamily{#3}\fontseries{#4}\fontshape{#5}%
  \selectfont}%
\fi\endgroup%
\begin{picture}(4074,834)(79,-253)
\put(1351,119){\makebox(0,0)[b]{\smash{{\SetFigFont{10}{12.0}{\rmdefault}{\mddefault}{\updefault}{\color[rgb]{0,0,0}$\equiv$}%
}}}}
\put(856,-61){\makebox(0,0)[b]{\smash{{\SetFigFont{10}{12.0}{\rmdefault}{\mddefault}{\updefault}{\color[rgb]{0,0,0}$\phi'$}%
}}}}
\put(316,299){\makebox(0,0)[b]{\smash{{\SetFigFont{10}{12.0}{\rmdefault}{\mddefault}{\updefault}{\color[rgb]{0,0,0}$\phi$}%
}}}}
\put(1846,119){\makebox(0,0)[b]{\smash{{\SetFigFont{10}{12.0}{\rmdefault}{\mddefault}{\updefault}{\color[rgb]{0,0,0}$\phi$}%
}}}}
\put(2386,119){\makebox(0,0)[b]{\smash{{\SetFigFont{10}{12.0}{\rmdefault}{\mddefault}{\updefault}{\color[rgb]{0,0,0}$\phi'$}%
}}}}
\put(3916,299){\makebox(0,0)[b]{\smash{{\SetFigFont{10}{12.0}{\rmdefault}{\mddefault}{\updefault}{\color[rgb]{0,0,0}$\phi'$}%
}}}}
\put(3376,-61){\makebox(0,0)[b]{\smash{{\SetFigFont{10}{12.0}{\rmdefault}{\mddefault}{\updefault}{\color[rgb]{0,0,0}$\phi$}%
}}}}
\put(2881,119){\makebox(0,0)[b]{\smash{{\SetFigFont{10}{12.0}{\rmdefault}{\mddefault}{\updefault}{\color[rgb]{0,0,0}$\equiv$}%
}}}}
\end{picture}%
\end{center}
\end{rem}

\begin{ex}\label{ex:circuits}
Let us consider the Petri net $(X,R)$ from example \ref{ex:reseau_de_petri}. Its associated $2$-polygraph is made of one $0$-cell $\ast$, three $1$-cells $x$, $y$ and $z$ and eight $2$-cells pictured as:
\begin{center}
\begin{picture}(0,0)%
\includegraphics{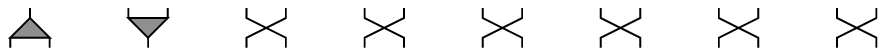}%
\end{picture}%
\setlength{\unitlength}{4144sp}%
\begingroup\makeatletter\ifx\SetFigFont\undefined%
\gdef\SetFigFont#1#2#3#4#5{%
  \reset@font\fontsize{#1}{#2pt}%
  \fontfamily{#3}\fontseries{#4}\fontshape{#5}%
  \selectfont}%
\fi\endgroup%
\begin{picture}(4191,754)(1685,-20)
\put(1981,254){\makebox(0,0)[b]{\smash{{\SetFigFont{10}{12.0}{\rmdefault}{\mddefault}{\updefault}{\color[rgb]{0,0,0}$z$}%
}}}}
\put(2341,614){\makebox(0,0)[b]{\smash{{\SetFigFont{10}{12.0}{\rmdefault}{\mddefault}{\updefault}{\color[rgb]{0,0,0}$y$}%
}}}}
\put(2521,614){\makebox(0,0)[b]{\smash{{\SetFigFont{10}{12.0}{\rmdefault}{\mddefault}{\updefault}{\color[rgb]{0,0,0}$y$}%
}}}}
\put(2431,254){\makebox(0,0)[b]{\smash{{\SetFigFont{10}{12.0}{\rmdefault}{\mddefault}{\updefault}{\color[rgb]{0,0,0}$z$}%
}}}}
\put(3061,254){\makebox(0,0)[b]{\smash{{\SetFigFont{10}{12.0}{\rmdefault}{\mddefault}{\updefault}{\color[rgb]{0,0,0}$x$}%
}}}}
\put(3061,614){\makebox(0,0)[b]{\smash{{\SetFigFont{10}{12.0}{\rmdefault}{\mddefault}{\updefault}{\color[rgb]{0,0,0}$y$}%
}}}}
\put(2881,614){\makebox(0,0)[b]{\smash{{\SetFigFont{10}{12.0}{\rmdefault}{\mddefault}{\updefault}{\color[rgb]{0,0,0}$x$}%
}}}}
\put(2881,254){\makebox(0,0)[b]{\smash{{\SetFigFont{10}{12.0}{\rmdefault}{\mddefault}{\updefault}{\color[rgb]{0,0,0}$y$}%
}}}}
\put(3421,254){\makebox(0,0)[b]{\smash{{\SetFigFont{10}{12.0}{\rmdefault}{\mddefault}{\updefault}{\color[rgb]{0,0,0}$x$}%
}}}}
\put(3601,614){\makebox(0,0)[b]{\smash{{\SetFigFont{10}{12.0}{\rmdefault}{\mddefault}{\updefault}{\color[rgb]{0,0,0}$x$}%
}}}}
\put(3421,614){\makebox(0,0)[b]{\smash{{\SetFigFont{10}{12.0}{\rmdefault}{\mddefault}{\updefault}{\color[rgb]{0,0,0}$y$}%
}}}}
\put(3601,254){\makebox(0,0)[b]{\smash{{\SetFigFont{10}{12.0}{\rmdefault}{\mddefault}{\updefault}{\color[rgb]{0,0,0}$y$}%
}}}}
\put(3961,614){\makebox(0,0)[b]{\smash{{\SetFigFont{10}{12.0}{\rmdefault}{\mddefault}{\updefault}{\color[rgb]{0,0,0}$x$}%
}}}}
\put(4141,254){\makebox(0,0)[b]{\smash{{\SetFigFont{10}{12.0}{\rmdefault}{\mddefault}{\updefault}{\color[rgb]{0,0,0}$x$}%
}}}}
\put(3961,254){\makebox(0,0)[b]{\smash{{\SetFigFont{10}{12.0}{\rmdefault}{\mddefault}{\updefault}{\color[rgb]{0,0,0}$z$}%
}}}}
\put(4141,614){\makebox(0,0)[b]{\smash{{\SetFigFont{10}{12.0}{\rmdefault}{\mddefault}{\updefault}{\color[rgb]{0,0,0}$z$}%
}}}}
\put(5041,614){\makebox(0,0)[b]{\smash{{\SetFigFont{10}{12.0}{\rmdefault}{\mddefault}{\updefault}{\color[rgb]{0,0,0}$y$}%
}}}}
\put(5221,254){\makebox(0,0)[b]{\smash{{\SetFigFont{10}{12.0}{\rmdefault}{\mddefault}{\updefault}{\color[rgb]{0,0,0}$y$}%
}}}}
\put(5041,254){\makebox(0,0)[b]{\smash{{\SetFigFont{10}{12.0}{\rmdefault}{\mddefault}{\updefault}{\color[rgb]{0,0,0}$z$}%
}}}}
\put(5221,614){\makebox(0,0)[b]{\smash{{\SetFigFont{10}{12.0}{\rmdefault}{\mddefault}{\updefault}{\color[rgb]{0,0,0}$z$}%
}}}}
\put(5581,254){\makebox(0,0)[b]{\smash{{\SetFigFont{10}{12.0}{\rmdefault}{\mddefault}{\updefault}{\color[rgb]{0,0,0}$y$}%
}}}}
\put(5761,614){\makebox(0,0)[b]{\smash{{\SetFigFont{10}{12.0}{\rmdefault}{\mddefault}{\updefault}{\color[rgb]{0,0,0}$y$}%
}}}}
\put(5581,614){\makebox(0,0)[b]{\smash{{\SetFigFont{10}{12.0}{\rmdefault}{\mddefault}{\updefault}{\color[rgb]{0,0,0}$z$}%
}}}}
\put(5761,254){\makebox(0,0)[b]{\smash{{\SetFigFont{10}{12.0}{\rmdefault}{\mddefault}{\updefault}{\color[rgb]{0,0,0}$z$}%
}}}}
\put(4681,614){\makebox(0,0)[b]{\smash{{\SetFigFont{10}{12.0}{\rmdefault}{\mddefault}{\updefault}{\color[rgb]{0,0,0}$x$}%
}}}}
\put(4501,254){\makebox(0,0)[b]{\smash{{\SetFigFont{10}{12.0}{\rmdefault}{\mddefault}{\updefault}{\color[rgb]{0,0,0}$x$}%
}}}}
\put(4501,614){\makebox(0,0)[b]{\smash{{\SetFigFont{10}{12.0}{\rmdefault}{\mddefault}{\updefault}{\color[rgb]{0,0,0}$z$}%
}}}}
\put(4681,254){\makebox(0,0)[b]{\smash{{\SetFigFont{10}{12.0}{\rmdefault}{\mddefault}{\updefault}{\color[rgb]{0,0,0}$z$}%
}}}}
\put(1891, 29){\makebox(0,0)[b]{\smash{{\SetFigFont{10}{12.0}{\rmdefault}{\mddefault}{\updefault}{\color[rgb]{0,0,0}$\alpha$}%
}}}}
\put(2431, 29){\makebox(0,0)[b]{\smash{{\SetFigFont{10}{12.0}{\rmdefault}{\mddefault}{\updefault}{\color[rgb]{0,0,0}$\beta$}%
}}}}
\put(2971, 29){\makebox(0,0)[b]{\smash{{\SetFigFont{10}{12.0}{\rmdefault}{\mddefault}{\updefault}{\color[rgb]{0,0,0}$\tau_{x,y}$}%
}}}}
\put(3511, 29){\makebox(0,0)[b]{\smash{{\SetFigFont{10}{12.0}{\rmdefault}{\mddefault}{\updefault}{\color[rgb]{0,0,0}$\tau_{y,x}$}%
}}}}
\put(4051, 29){\makebox(0,0)[b]{\smash{{\SetFigFont{10}{12.0}{\rmdefault}{\mddefault}{\updefault}{\color[rgb]{0,0,0}$\tau_{x,z}$}%
}}}}
\put(4591, 29){\makebox(0,0)[b]{\smash{{\SetFigFont{10}{12.0}{\rmdefault}{\mddefault}{\updefault}{\color[rgb]{0,0,0}$\tau_{z,x}$}%
}}}}
\put(5131, 29){\makebox(0,0)[b]{\smash{{\SetFigFont{10}{12.0}{\rmdefault}{\mddefault}{\updefault}{\color[rgb]{0,0,0}$\tau_{y,z}$}%
}}}}
\put(5671, 29){\makebox(0,0)[b]{\smash{{\SetFigFont{10}{12.0}{\rmdefault}{\mddefault}{\updefault}{\color[rgb]{0,0,0}$\tau_{z,y}$}%
}}}}
\put(1891,614){\makebox(0,0)[b]{\smash{{\SetFigFont{10}{12.0}{\rmdefault}{\mddefault}{\updefault}{\color[rgb]{0,0,0}$x$}%
}}}}
\put(1801,254){\makebox(0,0)[b]{\smash{{\SetFigFont{10}{12.0}{\rmdefault}{\mddefault}{\updefault}{\color[rgb]{0,0,0}$y$}%
}}}}
\end{picture}%
\end{center}

\noindent Then one considers the reduction graph from examples \ref{ex:graphe_de_reduction} and \ref{ex:mots_commutatifs}. As we have seen, all the paths in this diagram can be lifted to representatives in the free category $\mon{G(\Sigma^2(X,R)}$. These representatives are organized in a diagram such as the following one:
$$
\xymatrix
{
	x^2\tens y^2 
		\ar[rrr]^-{x^2\tens\beta} 
		\ar[d]_-{\alpha\tens x\tens y^2}
	&&& 
	x^2\tens z 
		\ar[d]^-{\alpha\tens x\tens z}
	\\
	y\tens z\tens x\tens y^2 
		\ar[rrr]^-{y\tens z\tens x\tens\beta}
		\ar[d]_-{y\tens z\tens\alpha\tens y^2}
	&&&
	y\tens z\tens x\tens z
		\ar[d]^-{y\tens z\tens\alpha\tens z}
  \\
	y\tens z\tens y\tens z\tens y^2 
		\ar[rrr]^-{y\tens z\tens y\tens z\tens\beta}
		\ar[d]_-{y\tens\tau_{z,y}\tens z\tens y^2}
	&&&
	y\tens z\tens y\tens z^2
		\ar[d]^-{y\tens\tau_{z,y}\tens z^2}
	\\
	y^2\tens z^2\tens y^2 
		\ar[rrr]^-{y^2\tens z^2\tens\beta}
		\ar[d]_-{\beta\tens z^2\tens y^2}
	&&&
	y^2\tens z^3
		\ar[d]^-{\beta\tens z^2\tens y^2}
	\\
	z^3\tens y^2 
		\ar[rrr]^-{z^3\tens\beta}
	&&&
	z^4.
}
$$

\noindent In this diagram, all parallel paths only differ by the order of application of the same $2$-cells in different parts of the same $1$-arrows: hence they are identified by the exchange relations, which means that they become equal in the free $2$-category generated by $\Sigma^2(X,R)$. For example, the $2$-arrow corresponding to any composite from $x^2\tens y^2$ to $z^4$ is written as $(\beta\tens z^3)\circ(y\tens\tau_{z,y}\tens z^2)\circ(\alpha\tens\alpha\tens\beta)$ and is pictured as the following more-readable circuit:
\begin{center}
\begin{picture}(0,0)%
\includegraphics{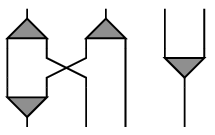}%
\end{picture}%
\setlength{\unitlength}{4144sp}%
\begingroup\makeatletter\ifx\SetFigFont\undefined%
\gdef\SetFigFont#1#2#3#4#5{%
  \reset@font\fontsize{#1}{#2pt}%
  \fontfamily{#3}\fontseries{#4}\fontshape{#5}%
  \selectfont}%
\fi\endgroup%
\begin{picture}(1049,970)(1317,-1001)
\put(1801,-151){\makebox(0,0)[b]{\smash{{\SetFigFont{10}{12.0}{\rmdefault}{\mddefault}{\updefault}{\color[rgb]{0,0,0}$x$}%
}}}}
\put(1441,-151){\makebox(0,0)[b]{\smash{{\SetFigFont{10}{12.0}{\rmdefault}{\mddefault}{\updefault}{\color[rgb]{0,0,0}$x$}%
}}}}
\put(2071,-151){\makebox(0,0)[b]{\smash{{\SetFigFont{10}{12.0}{\rmdefault}{\mddefault}{\updefault}{\color[rgb]{0,0,0}$y$}%
}}}}
\put(2251,-151){\makebox(0,0)[b]{\smash{{\SetFigFont{10}{12.0}{\rmdefault}{\mddefault}{\updefault}{\color[rgb]{0,0,0}$y$}%
}}}}
\put(1441,-961){\makebox(0,0)[b]{\smash{{\SetFigFont{10}{12.0}{\rmdefault}{\mddefault}{\updefault}{\color[rgb]{0,0,0}$z$}%
}}}}
\put(1711,-961){\makebox(0,0)[b]{\smash{{\SetFigFont{10}{12.0}{\rmdefault}{\mddefault}{\updefault}{\color[rgb]{0,0,0}$z$}%
}}}}
\put(1891,-961){\makebox(0,0)[b]{\smash{{\SetFigFont{10}{12.0}{\rmdefault}{\mddefault}{\updefault}{\color[rgb]{0,0,0}$z$}%
}}}}
\put(2161,-961){\makebox(0,0)[b]{\smash{{\SetFigFont{10}{12.0}{\rmdefault}{\mddefault}{\updefault}{\color[rgb]{0,0,0}$z$}%
}}}}
\end{picture}%
\end{center}
\end{ex}

\noindent From this example, it seems that the congruences $\equi{(X,R)}$ in $\mon{G(X,R)}$ and $\equi{01}$ in $\mon{G(\Sigma^2(X,R))}$ are linked in some way. For that, we  denote by $\Gr(X,R)$ the quotient category $\mon{G(X,R)}/\equi{(X,R)}$.

\begin{prop}\label{prop:2d-equiv}
Let $(X,R)$ be a Petri net. The functor $\pi:\mon{G(\Sigma^2(\Sigma,R))}\fl\mon{G(X,R)}$ induces a functor $\ol{\pi}:\mon{\Sigma^2(X,R)}\fl\Gr(X,R)$.
\end{prop}

\begin{dem}
We have to check that, whenever $f$ and $g$ are parallel arrows in $G(\Sigma^2(X,R))$ such that $f\equi{01} g$, we have $\pi(f)\equi{(X,R)}\pi(g)$. Let $u$, $v$, $w$ be $1$-arrows and $\alpha$, $\beta$ be $2$-cells in $\Sigma^2(X,R)$. Then, by definition of the functor~$\pi$, the following four equalities hold:
$$
\left\{
\begin{array}{r c l}
\pi(u\tens\alpha\tens v\tens s_1(\beta)\tens w) &\:=\: &(u+v+w+s_1(\beta))+\alpha, \\
\pi(u\tens\alpha\tens v\tens t_1(\beta)\tens w) &\:=\: &(u+v+w+t_1(\beta))+\alpha, \\
\pi(u\tens s_1(\alpha)\tens v\tens\beta\tens w) &\:=\: &(u+v+w+s_1(\alpha))+\beta, \\
\pi(u\tens t_1(\alpha)\tens v\tens\beta\tens w) &\:=\: &(u+v+w+t_1(\alpha))+\beta. 
\end{array}
\right.
$$

\noindent Thus, the functor $\pi$ satisfies:
$$
\pi\left(
\begin{array}{c}
u\tens\phi\tens(v\tens s_1(\psi)\tens w) \\
\circ \\
(u\tens t_1(\phi)\tens v)\tens\psi\tens w 
\end{array}
\right)
\qquad\equi{(X,R)}\qquad
\pi\left(
\begin{array}{c}
(u\tens s_1(\phi)\tens v)\tens\psi\tens w \\
\circ \\
u\tens\phi\tens(v\tens t_1(\psi)\tens w)
\end{array}
\right).
$$

\noindent Since $\pi$ is a functor, we get that $\pi(f)\equi{(X,R)}\pi(g)$ for any two parallel $f$ and $g$ such that $f\equi{01}g$. 

\findem\end{dem}

\noindent For the moment, we have seen that Petri nets can be translated as $2$-polygraphs $\Sigma=(\ast,\Sigma_1,\Sigma_2\amalg S_{\Sigma_1})$ where $\Sigma_1$ and $\Sigma_2$ are finite sets and where $S_X$ denotes the set of all $2$-cells $\tau_{x,y}$, with $x$ and $y$ distinct elements in~$X$. 

Conversely, given any $2$-polygraph of the form $\Sigma=(\ast,\Sigma_1,\Sigma_2\amalg S_{\Sigma_1})$ with $\Sigma_1$ and $\Sigma_2$ finite, one can build a Petri net with alphabet $\Sigma_1$ and rules given by the projection through $\pi:\mon{\Sigma_1}\fl[\Sigma_1]$ of the $2$-cells of $\Sigma_2$. Furthermore, it can be proved that the two transformations between Petri nets and $2$-polygraphs of this form are inverse to each other. 

Hence, we could state that Petri nets \emph{are} $2$-polygraphs of the form $\Sigma=(\ast,\Sigma_1,\Sigma_2\amalg S_{\Sigma_1})$. However, this would be quite excessive since there are much more $2$-arrows in $\mon{\Sigma}$ than rewriting paths in the corresponding Petri net. 

\begin{ex}
Once again, let us consider the Petri net from example \ref{ex:reseau_de_petri} and the path in $\Gr(X,R)$ given in examples \ref{ex:graphe_de_reduction} and \ref{ex:mots_commutatifs}. In example \ref{ex:circuits}, we have already seen a $2$-arrow of $\mon{\Sigma^2(X,R)}$ representing this reduction path. The following parallel $2$-arrows are also possible representatives for this path:
\begin{center}
\begin{picture}(0,0)%
\includegraphics{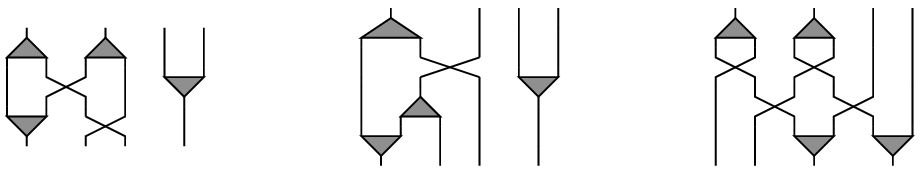}%
\end{picture}%
\setlength{\unitlength}{4144sp}%
\begingroup\makeatletter\ifx\SetFigFont\undefined%
\gdef\SetFigFont#1#2#3#4#5{%
  \reset@font\fontsize{#1}{#2pt}%
  \fontfamily{#3}\fontseries{#4}\fontshape{#5}%
  \selectfont}%
\fi\endgroup%
\begin{picture}(4288,1150)(-122,-281)
\put(2341,-241){\makebox(0,0)[b]{\smash{{\SetFigFont{10}{12.0}{\rmdefault}{\mddefault}{\updefault}{\color[rgb]{0,0,0}$z$}%
}}}}
\put(3241,749){\makebox(0,0)[b]{\smash{{\SetFigFont{10}{12.0}{\rmdefault}{\mddefault}{\updefault}{\color[rgb]{0,0,0}$x$}%
}}}}
\put(3601,749){\makebox(0,0)[b]{\smash{{\SetFigFont{10}{12.0}{\rmdefault}{\mddefault}{\updefault}{\color[rgb]{0,0,0}$x$}%
}}}}
\put(3871,749){\makebox(0,0)[b]{\smash{{\SetFigFont{10}{12.0}{\rmdefault}{\mddefault}{\updefault}{\color[rgb]{0,0,0}$y$}%
}}}}
\put(4051,749){\makebox(0,0)[b]{\smash{{\SetFigFont{10}{12.0}{\rmdefault}{\mddefault}{\updefault}{\color[rgb]{0,0,0}$y$}%
}}}}
\put(3151,-241){\makebox(0,0)[b]{\smash{{\SetFigFont{10}{12.0}{\rmdefault}{\mddefault}{\updefault}{\color[rgb]{0,0,0}$z$}%
}}}}
\put(3331,-241){\makebox(0,0)[b]{\smash{{\SetFigFont{10}{12.0}{\rmdefault}{\mddefault}{\updefault}{\color[rgb]{0,0,0}$z$}%
}}}}
\put(3601,-241){\makebox(0,0)[b]{\smash{{\SetFigFont{10}{12.0}{\rmdefault}{\mddefault}{\updefault}{\color[rgb]{0,0,0}$z$}%
}}}}
\put(3961,-241){\makebox(0,0)[b]{\smash{{\SetFigFont{10}{12.0}{\rmdefault}{\mddefault}{\updefault}{\color[rgb]{0,0,0}$z$}%
}}}}
\put(  1,659){\makebox(0,0)[b]{\smash{{\SetFigFont{10}{12.0}{\rmdefault}{\mddefault}{\updefault}{\color[rgb]{0,0,0}$x$}%
}}}}
\put(631,659){\makebox(0,0)[b]{\smash{{\SetFigFont{10}{12.0}{\rmdefault}{\mddefault}{\updefault}{\color[rgb]{0,0,0}$y$}%
}}}}
\put(811,659){\makebox(0,0)[b]{\smash{{\SetFigFont{10}{12.0}{\rmdefault}{\mddefault}{\updefault}{\color[rgb]{0,0,0}$y$}%
}}}}
\put(  1,-151){\makebox(0,0)[b]{\smash{{\SetFigFont{10}{12.0}{\rmdefault}{\mddefault}{\updefault}{\color[rgb]{0,0,0}$z$}%
}}}}
\put(271,-151){\makebox(0,0)[b]{\smash{{\SetFigFont{10}{12.0}{\rmdefault}{\mddefault}{\updefault}{\color[rgb]{0,0,0}$z$}%
}}}}
\put(451,-151){\makebox(0,0)[b]{\smash{{\SetFigFont{10}{12.0}{\rmdefault}{\mddefault}{\updefault}{\color[rgb]{0,0,0}$z$}%
}}}}
\put(721,-151){\makebox(0,0)[b]{\smash{{\SetFigFont{10}{12.0}{\rmdefault}{\mddefault}{\updefault}{\color[rgb]{0,0,0}$z$}%
}}}}
\put(361,659){\makebox(0,0)[b]{\smash{{\SetFigFont{10}{12.0}{\rmdefault}{\mddefault}{\updefault}{\color[rgb]{0,0,0}$x$}%
}}}}
\put(1666,749){\makebox(0,0)[b]{\smash{{\SetFigFont{10}{12.0}{\rmdefault}{\mddefault}{\updefault}{\color[rgb]{0,0,0}$x$}%
}}}}
\put(2071,749){\makebox(0,0)[b]{\smash{{\SetFigFont{10}{12.0}{\rmdefault}{\mddefault}{\updefault}{\color[rgb]{0,0,0}$x$}%
}}}}
\put(2251,749){\makebox(0,0)[b]{\smash{{\SetFigFont{10}{12.0}{\rmdefault}{\mddefault}{\updefault}{\color[rgb]{0,0,0}$y$}%
}}}}
\put(2431,749){\makebox(0,0)[b]{\smash{{\SetFigFont{10}{12.0}{\rmdefault}{\mddefault}{\updefault}{\color[rgb]{0,0,0}$y$}%
}}}}
\put(1621,-241){\makebox(0,0)[b]{\smash{{\SetFigFont{10}{12.0}{\rmdefault}{\mddefault}{\updefault}{\color[rgb]{0,0,0}$z$}%
}}}}
\put(1891,-241){\makebox(0,0)[b]{\smash{{\SetFigFont{10}{12.0}{\rmdefault}{\mddefault}{\updefault}{\color[rgb]{0,0,0}$z$}%
}}}}
\put(2071,-241){\makebox(0,0)[b]{\smash{{\SetFigFont{10}{12.0}{\rmdefault}{\mddefault}{\updefault}{\color[rgb]{0,0,0}$z$}%
}}}}
\end{picture}%
\end{center}
\end{ex}

\noindent Hence, even if there is a correspondance between Petri nets and $2$-polygraphs $\Sigma=(\ast,\Sigma_1,\Sigma_2\amalg S_{\Sigma_1})$, both objects do not naturally generate the same reduction graphs since $\mon{\Sigma^2(X,R)}$ is bigger than $\Gr(X,R)$. There are many possible solutions to this problem. One possibility is to add relations between parallel $2$-arrows of $\Sigma$ that represent the same path in the Petri net reduction graph: we are going to sketch such a study in the rest of this section. Another really different solution is studied in section~\ref{sec:3d}, where we use the fact that commutative monoids correspond to a special class of $2$-polygraphs.

For the moment, let us consider a $2$-polygraph $\Sigma=(\ast,\Sigma_1,\Sigma_2\amalg S_{\Sigma_1})$, but where $S_{\Sigma_1}$ now also contains explicit permutations $\tau_{x,x}$ for every $1$-cell~$x$ in~$\Sigma_1$. This extension does not change the properties studied so far if we extend the functor $\pi$ with $\pi(\tau_{x,x})=\id_{x+x}$. We denote by $(\Sigma_1,\Sigma_2)$ the corresponding Petri net. The following result gives a family of relations for some parallel $2$-arrows corresponding to the same Petri net reduction. Its proof is straightforward and uses the facts that $\ol{\pi}$ is a functor and maps each $\tau_{x,y}$ onto an identity.

\begin{lem}\label{lem:relations}
The functor $\ol{\pi}$ is compatible with the congruence $\equiv$ generated by the following relations, given for all $1$-cells~$x$,~$y$ and~$z$ and every $2$-cell $\alpha$:
\begin{center}
\begin{picture}(0,0)%
\includegraphics{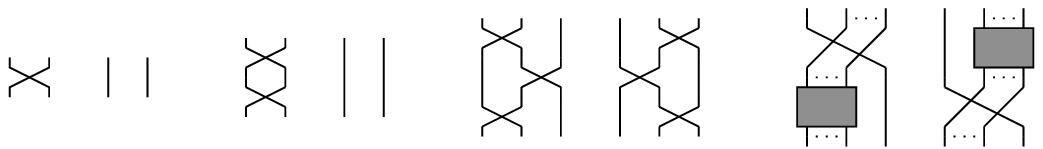}%
\end{picture}%
\setlength{\unitlength}{4144sp}%
\begingroup\makeatletter\ifx\SetFigFont\undefined%
\gdef\SetFigFont#1#2#3#4#5{%
  \reset@font\fontsize{#1}{#2pt}%
  \fontfamily{#3}\fontseries{#4}\fontshape{#5}%
  \selectfont}%
\fi\endgroup%
\begin{picture}(4815,795)(-1022,152)
\put(-584,434){\makebox(0,0)[b]{\smash{{\SetFigFont{10}{12.0}{\rmdefault}{\mddefault}{\updefault}{\color[rgb]{0,0,0}$\equiv$}%
}}}}
\put(2836,299){\makebox(0,0)[b]{\smash{{\SetFigFont{10}{12.0}{\rmdefault}{\mddefault}{\updefault}{\color[rgb]{0,0,0}$\alpha$}%
}}}}
\put(3646,569){\makebox(0,0)[b]{\smash{{\SetFigFont{10}{12.0}{\rmdefault}{\mddefault}{\updefault}{\color[rgb]{0,0,0}$\alpha$}%
}}}}
\put(181,704){\makebox(0,0)[b]{\smash{{\SetFigFont{10}{12.0}{\rmdefault}{\mddefault}{\updefault}{\color[rgb]{0,0,0}$x$}%
}}}}
\put(361,704){\makebox(0,0)[b]{\smash{{\SetFigFont{10}{12.0}{\rmdefault}{\mddefault}{\updefault}{\color[rgb]{0,0,0}$y$}%
}}}}
\put(631,704){\makebox(0,0)[b]{\smash{{\SetFigFont{10}{12.0}{\rmdefault}{\mddefault}{\updefault}{\color[rgb]{0,0,0}$x$}%
}}}}
\put(811,704){\makebox(0,0)[b]{\smash{{\SetFigFont{10}{12.0}{\rmdefault}{\mddefault}{\updefault}{\color[rgb]{0,0,0}$y$}%
}}}}
\put(1261,794){\makebox(0,0)[b]{\smash{{\SetFigFont{10}{12.0}{\rmdefault}{\mddefault}{\updefault}{\color[rgb]{0,0,0}$x$}%
}}}}
\put(1441,794){\makebox(0,0)[b]{\smash{{\SetFigFont{10}{12.0}{\rmdefault}{\mddefault}{\updefault}{\color[rgb]{0,0,0}$y$}%
}}}}
\put(1621,794){\makebox(0,0)[b]{\smash{{\SetFigFont{10}{12.0}{\rmdefault}{\mddefault}{\updefault}{\color[rgb]{0,0,0}$z$}%
}}}}
\put(1891,794){\makebox(0,0)[b]{\smash{{\SetFigFont{10}{12.0}{\rmdefault}{\mddefault}{\updefault}{\color[rgb]{0,0,0}$x$}%
}}}}
\put(2071,794){\makebox(0,0)[b]{\smash{{\SetFigFont{10}{12.0}{\rmdefault}{\mddefault}{\updefault}{\color[rgb]{0,0,0}$y$}%
}}}}
\put(2251,794){\makebox(0,0)[b]{\smash{{\SetFigFont{10}{12.0}{\rmdefault}{\mddefault}{\updefault}{\color[rgb]{0,0,0}$z$}%
}}}}
\put(2746,839){\makebox(0,0)[b]{\smash{{\SetFigFont{10}{12.0}{\rmdefault}{\mddefault}{\updefault}{\color[rgb]{0,0,0}$x$}%
}}}}
\put(3376,839){\makebox(0,0)[b]{\smash{{\SetFigFont{10}{12.0}{\rmdefault}{\mddefault}{\updefault}{\color[rgb]{0,0,0}$x$}%
}}}}
\put(496,434){\makebox(0,0)[b]{\smash{{\SetFigFont{10}{12.0}{\rmdefault}{\mddefault}{\updefault}{\color[rgb]{0,0,0}$\equiv$}%
}}}}
\put(1756,434){\makebox(0,0)[b]{\smash{{\SetFigFont{10}{12.0}{\rmdefault}{\mddefault}{\updefault}{\color[rgb]{0,0,0}$\equiv$}%
}}}}
\put(3241,434){\makebox(0,0)[b]{\smash{{\SetFigFont{10}{12.0}{\rmdefault}{\mddefault}{\updefault}{\color[rgb]{0,0,0}$\equiv$}%
}}}}
\put(-899,614){\makebox(0,0)[b]{\smash{{\SetFigFont{10}{12.0}{\rmdefault}{\mddefault}{\updefault}{\color[rgb]{0,0,0}$x$}%
}}}}
\put(-719,614){\makebox(0,0)[b]{\smash{{\SetFigFont{10}{12.0}{\rmdefault}{\mddefault}{\updefault}{\color[rgb]{0,0,0}$x$}%
}}}}
\put(-449,614){\makebox(0,0)[b]{\smash{{\SetFigFont{10}{12.0}{\rmdefault}{\mddefault}{\updefault}{\color[rgb]{0,0,0}$x$}%
}}}}
\put(-269,614){\makebox(0,0)[b]{\smash{{\SetFigFont{10}{12.0}{\rmdefault}{\mddefault}{\updefault}{\color[rgb]{0,0,0}$x$}%
}}}}
\end{picture}%
\end{center}

\noindent where the generalized explicit permutations used in the third relation are defined inductively from the permutation $2$-cells in a graphically intuitive way.
\end{lem}

\begin{rem}
The first relation states that, in a given marking of a Petri net, two tokens in the same place are totally indiscernible: for example, one cannot tell if a given transition has consumed one given token or another one in the same place.
\end{rem}

\begin{rem}
We conjecture that the congruence $\equiv$ also satisfies the converse property: if $f$ and $g$ are two parallel $2$-arrows in $\mon{\Sigma}$ such that $\ol{\pi}(f)=\ol{\pi}(g)$, then $f\equiv g$. However, we do not yet have a proof of this fact.
\end{rem}

\noindent So far, we have a set of equations relating $2$-arrows we wish to identify. However this raises a $2$-\emph{dimensional word problem} [Burroni 1993]: given two parallel $2$-arrows in $\mon{\Sigma}$, are they equal \emph{modulo} the congruence $\equiv$ or not? One way to build a decision procedure for such a problem is to follow the methodology developped in [Lafont 2003] and [Guiraud 2004] and build a convergent $3$-\emph{polygraph} equivalent to the given equational presentation. 

\begin{rem}
Here, we do not recall basic notions about rewriting: they can be found in [Baader Nipkow 1998] for example. Let us say that, for this section, a $3$-polygraph is specified by a $2$-polygraph equipped with rewriting rules between parallel $2$-arrows. These rules are in fact $3$-cells, but we postpone all definitions until section~\ref{sec:3d} since we only need the intuition of it being a "circuit rewriting system" here.
\end{rem}

\noindent We would like to craft a convergent $3$-polygraph for the congruence $\equiv$ on the $2$-category $\mon{\Sigma}$. However, the fact that $2$-cells may have several inputs \emph{and} several outputs at the same time makes the rewriting study much different than in the already-encountered cases. We give here a possible starting point for future work.

\begin{rem}
For this introduction, we limit ourselves on several points:
\begin{enumerate}
\item[-] First of all, we only consider the congruence $\equi{0}$ generated by the last third families: we remove the relations $\tau_{x,x}\equiv\id_{x\tens x}$ since we still do not know how to handle them. This must be seen as a first step towards the study of $\equiv$.
\item[-] The second limitation is that we assume that $\Sigma_2$ does not contain any $2$-cell with an empty output: the corresponding Petri net cannot have any transition that do not produce any token.
\item[-] Finally, we suppose that every $2$-cell in $\Sigma_2$ with an empty input has only one output. This is not a real limitation since, in a Petri net, we can replace a transition $\alpha:\ast\fl y_1+\dots+y_n$ by two transitions $\ast\fl z$ and $z\fl y_1+\dots+y_n$, with $z$ a new place. The Petri net one gets fully simulates the original one.
\end{enumerate}
\end{rem}

\noindent The idea is the following one: instead of giving an answer to the question $f\equiv g$ directly in $\mon{\Sigma}$, we translate $2$-arrows of $\Sigma$ into a $3$-polygraph in which we know a decision procedure and such that the translation preserves the congruence $\equiv$.

\begin{ntt}
We denote by $\ol{\Sigma}$ the $2$-polygraph with one cell in dimension $0$, with $\Sigma_1$ as its set of $1$-cells and with the following families of $2$-cells:
\begin{center}
\begin{picture}(0,0)%
\includegraphics{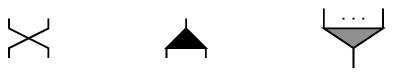}%
\end{picture}%
\setlength{\unitlength}{4144sp}%
\begingroup\makeatletter\ifx\SetFigFont\undefined%
\gdef\SetFigFont#1#2#3#4#5{%
  \reset@font\fontsize{#1}{#2pt}%
  \fontfamily{#3}\fontseries{#4}\fontshape{#5}%
  \selectfont}%
\fi\endgroup%
\begin{picture}(2033,721)(-32,295)
\put(1666,344){\makebox(0,0)[b]{\smash{{\SetFigFont{10}{12.0}{\rmdefault}{\mddefault}{\updefault}{\color[rgb]{0,0,0}$y_i$}%
}}}}
\put( 91,839){\makebox(0,0)[b]{\smash{{\SetFigFont{10}{12.0}{\rmdefault}{\mddefault}{\updefault}{\color[rgb]{0,0,0}$x$}%
}}}}
\put(271,839){\makebox(0,0)[b]{\smash{{\SetFigFont{10}{12.0}{\rmdefault}{\mddefault}{\updefault}{\color[rgb]{0,0,0}$y$}%
}}}}
\put( 91,389){\makebox(0,0)[b]{\smash{{\SetFigFont{10}{12.0}{\rmdefault}{\mddefault}{\updefault}{\color[rgb]{0,0,0}$y$}%
}}}}
\put(271,389){\makebox(0,0)[b]{\smash{{\SetFigFont{10}{12.0}{\rmdefault}{\mddefault}{\updefault}{\color[rgb]{0,0,0}$x$}%
}}}}
\put(901,839){\makebox(0,0)[b]{\smash{{\SetFigFont{10}{12.0}{\rmdefault}{\mddefault}{\updefault}{\color[rgb]{0,0,0}$x$}%
}}}}
\put(811,389){\makebox(0,0)[b]{\smash{{\SetFigFont{10}{12.0}{\rmdefault}{\mddefault}{\updefault}{\color[rgb]{0,0,0}$x$}%
}}}}
\put(991,389){\makebox(0,0)[b]{\smash{{\SetFigFont{10}{12.0}{\rmdefault}{\mddefault}{\updefault}{\color[rgb]{0,0,0}$x$}%
}}}}
\put(1531,884){\makebox(0,0)[b]{\smash{{\SetFigFont{10}{12.0}{\rmdefault}{\mddefault}{\updefault}{\color[rgb]{0,0,0}$x_1$}%
}}}}
\put(1801,884){\makebox(0,0)[b]{\smash{{\SetFigFont{10}{12.0}{\rmdefault}{\mddefault}{\updefault}{\color[rgb]{0,0,0}$x_n$}%
}}}}
\end{picture}%
\end{center}

\noindent The first family $(\tau_{x,y})$ is indexed by every possible $1$-cells $x$ and $y$; the second family $(\delta_x)$ by every $1$-cell $x$; the last one $(\alpha_i)$ by every $2$-cell $\alpha:x_1\dots x_m\fl y_1\dots y_n$ and every $i$ in $\dens{1}{n}$.

On top of the $2$-category $\mon{\ol{\Sigma}}$, we denote by $R$ the family made of the following $3$-cells, given for all possible coloration of the wires by $1$-cells:
\begin{center}
\begin{picture}(0,0)%
\includegraphics{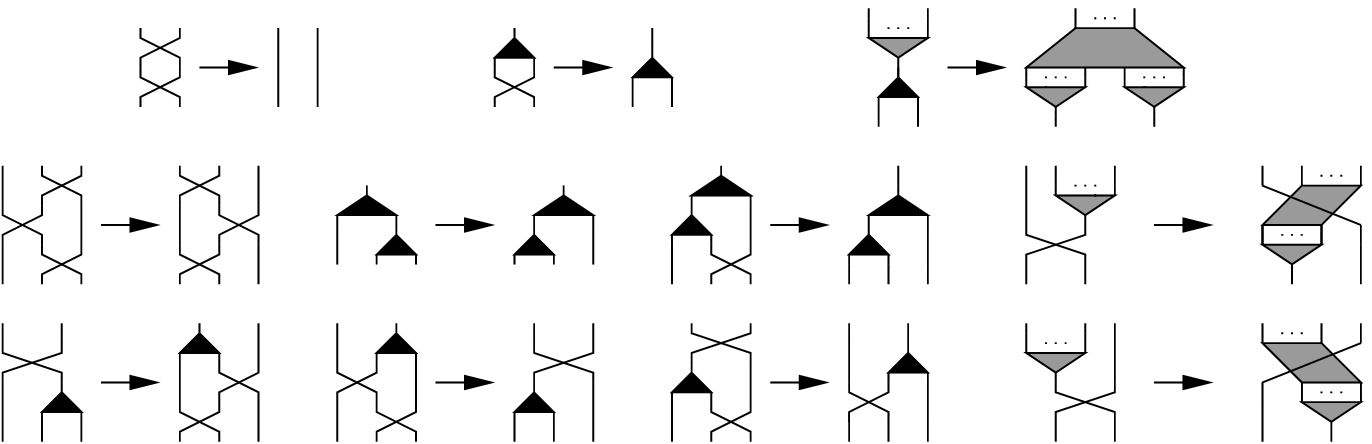}%
\end{picture}%
\setlength{\unitlength}{4144sp}%
\begingroup\makeatletter\ifx\SetFigFont\undefined%
\gdef\SetFigFont#1#2#3#4#5{%
  \reset@font\fontsize{#1}{#2pt}%
  \fontfamily{#3}\fontseries{#4}\fontshape{#5}%
  \selectfont}%
\fi\endgroup%
\begin{picture}(6312,2004)(79,-1243)
\put(6391,-1141){\makebox(0,0)[b]{\smash{{\SetFigFont{10}{12.0}{\familydefault}{\mddefault}{\updefault}{\color[rgb]{0,0,0}$\alpha_i$}%
}}}}
\put(5761,-511){\makebox(0,0)[b]{\smash{{\SetFigFont{10}{12.0}{\familydefault}{\mddefault}{\updefault}{\color[rgb]{0,0,0}$\alpha_i$}%
}}}}
\put(5266,-151){\makebox(0,0)[b]{\smash{{\SetFigFont{10}{12.0}{\familydefault}{\mddefault}{\updefault}{\color[rgb]{0,0,0}$\alpha_i$}%
}}}}
\put(5491,209){\makebox(0,0)[b]{\smash{{\SetFigFont{10}{12.0}{\familydefault}{\mddefault}{\updefault}{\color[rgb]{0,0,0}$\alpha_i$}%
}}}}
\put(3961,524){\makebox(0,0)[b]{\smash{{\SetFigFont{10}{12.0}{\familydefault}{\mddefault}{\updefault}{\color[rgb]{0,0,0}$\alpha_i$}%
}}}}
\put(4771,209){\makebox(0,0)[b]{\smash{{\SetFigFont{10}{12.0}{\familydefault}{\mddefault}{\updefault}{\color[rgb]{0,0,0}$\alpha_i$}%
}}}}
\put(4681,-916){\makebox(0,0)[b]{\smash{{\SetFigFont{10}{12.0}{\familydefault}{\mddefault}{\updefault}{\color[rgb]{0,0,0}$\alpha_i$}%
}}}}
\end{picture}%
\end{center}

\noindent The generalized duplication in the topmost-rightmost family is inductively built from local duplications and local permutations in a inductive way described in [Guiraud 2004] for example. We denote by $\equi{R}$ the congruence relation generated by $R$ on parallel $2$-arrows of $\mon{\ol{\Sigma}}$. 
\end{ntt}

\noindent Following the same method as the one presented in [Guiraud 2004] and using the coloration technique sketched in [Guiraud 2005], one proves that the $3$-polygraph $\ol{\Sigma}$ is convergent. Hence, given parallel $2$-arrows $f$ and $g$ in $\mon{\ol{\Sigma}}$, one can decide whether $f\equi{R}g$ holds or not.

Furthermore, we conjecture here that it is possible to define a $2$-functor $\Phi:\mon{\Sigma}\fl\mon{\ol{\Sigma}}$ such that $f\equi{0}g$ holds if and only if $\Phi(f)\equi{R}\Phi(g)$ holds. Here we define a $2$-functor $\Phi$ which is a good candidate for this rôle and check the easy part of the claim.

\begin{ntt}
We define a $2$-functor $\Phi:\mon{\Sigma}\fl\mon{\ol{\Sigma}}$ by giving its values on the cells of $\Sigma$:
\begin{enumerate}
\item[0.] It sends the only $0$-cell of $\Sigma$ onto itself.
\item[1.] It sends each $1$-cell $x$ of $\Sigma$ onto itself.
\item[2.] It sends each $\tau_{x,y}$ onto itself and, for every $2$-cell $\alpha:x_1\dots x_m\fl y_1\dots y_n$ in $\Sigma_2$ with $n\geq 2$, we define:
$$
\Phi(\alpha) \:=\: (\alpha_1\tens\dots\tens\alpha_n)\circ\delta^n_{x_1\dots x_m}
$$

\noindent where $\delta^n_{x_1\dots x_m}$ is the only generalized duplication from $x_1\dots x_m$ to $(x_1\dots x_m)^n$ that is in normal form with respect to~$R$. 
\end{enumerate}
\end{ntt}

\noindent Then we have:

\begin{prop}
The congruence $\Phi(\equi{0})$ is included into $\equi{R}$. 
\end{prop}

\begin{dem}
We check that, for every relation $f\equiv g$ defining $\equi{0}$, we have $\Phi(f)\equi{R}\Phi(g)$. This is immediate for the two relations that only involve local permutations. And for the third family of equations:
\begin{center}
\begin{picture}(0,0)%
\includegraphics{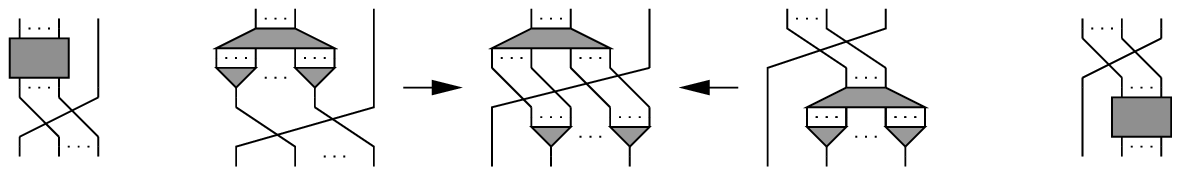}%
\end{picture}%
\setlength{\unitlength}{4144sp}%
\begingroup\makeatletter\ifx\SetFigFont\undefined%
\gdef\SetFigFont#1#2#3#4#5{%
  \reset@font\fontsize{#1}{#2pt}%
  \fontfamily{#3}\fontseries{#4}\fontshape{#5}%
  \selectfont}%
\fi\endgroup%
\begin{picture}(5867,885)(721,-118)
\put(6481,209){\makebox(0,0)[b]{\smash{{\SetFigFont{10}{12.0}{\rmdefault}{\mddefault}{\updefault}{\color[rgb]{0,0,0}$)$}%
}}}}
\put(1171,344){\makebox(0,0)[b]{\smash{{\SetFigFont{10}{12.0}{\rmdefault}{\mddefault}{\updefault}{\color[rgb]{0,0,0}$\alpha$}%
}}}}
\put(1441,614){\makebox(0,0)[b]{\smash{{\SetFigFont{10}{12.0}{\rmdefault}{\mddefault}{\updefault}{\color[rgb]{0,0,0}$x$}%
}}}}
\put(6211, 74){\makebox(0,0)[b]{\smash{{\SetFigFont{10}{12.0}{\rmdefault}{\mddefault}{\updefault}{\color[rgb]{0,0,0}$\alpha$}%
}}}}
\put(6301,614){\makebox(0,0)[b]{\smash{{\SetFigFont{10}{12.0}{\rmdefault}{\mddefault}{\updefault}{\color[rgb]{0,0,0}$x$}%
}}}}
\put(2701,659){\makebox(0,0)[b]{\smash{{\SetFigFont{10}{12.0}{\rmdefault}{\mddefault}{\updefault}{\color[rgb]{0,0,0}$x$}%
}}}}
\put(3961,659){\makebox(0,0)[b]{\smash{{\SetFigFont{10}{12.0}{\rmdefault}{\mddefault}{\updefault}{\color[rgb]{0,0,0}$x$}%
}}}}
\put(5041,659){\makebox(0,0)[b]{\smash{{\SetFigFont{10}{12.0}{\rmdefault}{\mddefault}{\updefault}{\color[rgb]{0,0,0}$x$}%
}}}}
\put(1801,209){\makebox(0,0)[b]{\smash{{\SetFigFont{10}{12.0}{\rmdefault}{\mddefault}{\updefault}{\color[rgb]{0,0,0}$=$}%
}}}}
\put(1621,209){\makebox(0,0)[b]{\smash{{\SetFigFont{10}{12.0}{\rmdefault}{\mddefault}{\updefault}{\color[rgb]{0,0,0}$)$}%
}}}}
\put(901,209){\makebox(0,0)[b]{\smash{{\SetFigFont{10}{12.0}{\rmdefault}{\mddefault}{\updefault}{\color[rgb]{0,0,0}$($}%
}}}}
\put(721,209){\makebox(0,0)[b]{\smash{{\SetFigFont{10}{12.0}{\rmdefault}{\mddefault}{\updefault}{\color[rgb]{0,0,0}$\Phi$}%
}}}}
\put(5401,209){\makebox(0,0)[b]{\smash{{\SetFigFont{10}{12.0}{\rmdefault}{\mddefault}{\updefault}{\color[rgb]{0,0,0}$=$}%
}}}}
\put(5581,209){\makebox(0,0)[b]{\smash{{\SetFigFont{10}{12.0}{\rmdefault}{\mddefault}{\updefault}{\color[rgb]{0,0,0}$\Phi$}%
}}}}
\put(5761,209){\makebox(0,0)[b]{\smash{{\SetFigFont{10}{12.0}{\rmdefault}{\mddefault}{\updefault}{\color[rgb]{0,0,0}$($}%
}}}}
\end{picture}%
\end{center}
\findem\end{dem}

\begin{ex}
Let us consider the Petri net from example \ref{ex:reseau_de_petri}. Its associated $2$-polygraph $\ol{\Sigma}$ has the following $2$-cells, beside the nine explicit permutations $(\tau_{\xi,\xi'})_{\xi,\xi'\in\ens{x,y,z}}$:
\begin{center}
\begin{picture}(0,0)%
\includegraphics{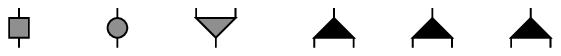}%
\end{picture}%
\setlength{\unitlength}{4144sp}%
\begingroup\makeatletter\ifx\SetFigFont\undefined%
\gdef\SetFigFont#1#2#3#4#5{%
  \reset@font\fontsize{#1}{#2pt}%
  \fontfamily{#3}\fontseries{#4}\fontshape{#5}%
  \selectfont}%
\fi\endgroup%
\begin{picture}(2669,754)(1407,-20)
\put(3871, 29){\makebox(0,0)[b]{\smash{{\SetFigFont{10}{12.0}{\rmdefault}{\mddefault}{\updefault}{\color[rgb]{0,0,0}$\delta_z$}%
}}}}
\put(2341,614){\makebox(0,0)[b]{\smash{{\SetFigFont{10}{12.0}{\rmdefault}{\mddefault}{\updefault}{\color[rgb]{0,0,0}$y$}%
}}}}
\put(2521,614){\makebox(0,0)[b]{\smash{{\SetFigFont{10}{12.0}{\rmdefault}{\mddefault}{\updefault}{\color[rgb]{0,0,0}$y$}%
}}}}
\put(2431,254){\makebox(0,0)[b]{\smash{{\SetFigFont{10}{12.0}{\rmdefault}{\mddefault}{\updefault}{\color[rgb]{0,0,0}$z$}%
}}}}
\put(2431, 29){\makebox(0,0)[b]{\smash{{\SetFigFont{10}{12.0}{\rmdefault}{\mddefault}{\updefault}{\color[rgb]{0,0,0}$\beta$}%
}}}}
\put(1981,254){\makebox(0,0)[b]{\smash{{\SetFigFont{10}{12.0}{\rmdefault}{\mddefault}{\updefault}{\color[rgb]{0,0,0}$z$}%
}}}}
\put(1531,254){\makebox(0,0)[b]{\smash{{\SetFigFont{10}{12.0}{\rmdefault}{\mddefault}{\updefault}{\color[rgb]{0,0,0}$y$}%
}}}}
\put(1531,614){\makebox(0,0)[b]{\smash{{\SetFigFont{10}{12.0}{\rmdefault}{\mddefault}{\updefault}{\color[rgb]{0,0,0}$x$}%
}}}}
\put(1981,614){\makebox(0,0)[b]{\smash{{\SetFigFont{10}{12.0}{\rmdefault}{\mddefault}{\updefault}{\color[rgb]{0,0,0}$x$}%
}}}}
\put(1531, 29){\makebox(0,0)[b]{\smash{{\SetFigFont{10}{12.0}{\rmdefault}{\mddefault}{\updefault}{\color[rgb]{0,0,0}$\alpha_1$}%
}}}}
\put(1981, 29){\makebox(0,0)[b]{\smash{{\SetFigFont{10}{12.0}{\rmdefault}{\mddefault}{\updefault}{\color[rgb]{0,0,0}$\alpha_2$}%
}}}}
\put(2971,614){\makebox(0,0)[b]{\smash{{\SetFigFont{10}{12.0}{\rmdefault}{\mddefault}{\updefault}{\color[rgb]{0,0,0}$x$}%
}}}}
\put(3421,614){\makebox(0,0)[b]{\smash{{\SetFigFont{10}{12.0}{\rmdefault}{\mddefault}{\updefault}{\color[rgb]{0,0,0}$y$}%
}}}}
\put(3871,614){\makebox(0,0)[b]{\smash{{\SetFigFont{10}{12.0}{\rmdefault}{\mddefault}{\updefault}{\color[rgb]{0,0,0}$z$}%
}}}}
\put(2881,254){\makebox(0,0)[b]{\smash{{\SetFigFont{10}{12.0}{\rmdefault}{\mddefault}{\updefault}{\color[rgb]{0,0,0}$x$}%
}}}}
\put(3061,254){\makebox(0,0)[b]{\smash{{\SetFigFont{10}{12.0}{\rmdefault}{\mddefault}{\updefault}{\color[rgb]{0,0,0}$x$}%
}}}}
\put(3331,254){\makebox(0,0)[b]{\smash{{\SetFigFont{10}{12.0}{\rmdefault}{\mddefault}{\updefault}{\color[rgb]{0,0,0}$y$}%
}}}}
\put(3511,254){\makebox(0,0)[b]{\smash{{\SetFigFont{10}{12.0}{\rmdefault}{\mddefault}{\updefault}{\color[rgb]{0,0,0}$y$}%
}}}}
\put(3781,254){\makebox(0,0)[b]{\smash{{\SetFigFont{10}{12.0}{\rmdefault}{\mddefault}{\updefault}{\color[rgb]{0,0,0}$z$}%
}}}}
\put(3961,254){\makebox(0,0)[b]{\smash{{\SetFigFont{10}{12.0}{\rmdefault}{\mddefault}{\updefault}{\color[rgb]{0,0,0}$z$}%
}}}}
\put(2971, 29){\makebox(0,0)[b]{\smash{{\SetFigFont{10}{12.0}{\rmdefault}{\mddefault}{\updefault}{\color[rgb]{0,0,0}$\delta_x$}%
}}}}
\put(3421, 29){\makebox(0,0)[b]{\smash{{\SetFigFont{10}{12.0}{\rmdefault}{\mddefault}{\updefault}{\color[rgb]{0,0,0}$\delta_y$}%
}}}}
\end{picture}%
\end{center}

\noindent Once translated into $\ol{\Sigma}$, the four representative we have seen of the Petri net reduction of example~\ref{ex:graphe_de_reduction} have the following respective normal forms:
\begin{center}
\begin{picture}(0,0)%
\includegraphics{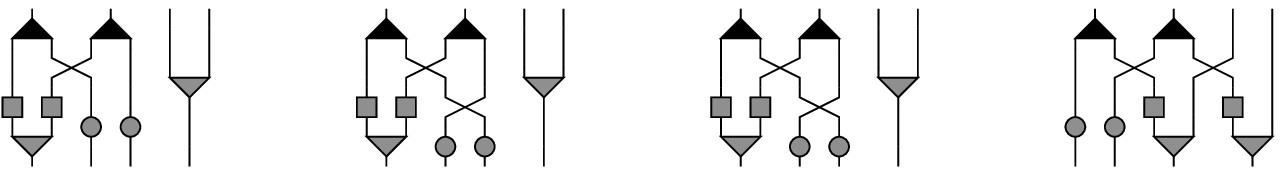}%
\end{picture}%
\setlength{\unitlength}{4144sp}%
\begingroup\makeatletter\ifx\SetFigFont\undefined%
\gdef\SetFigFont#1#2#3#4#5{%
  \reset@font\fontsize{#1}{#2pt}%
  \fontfamily{#3}\fontseries{#4}\fontshape{#5}%
  \selectfont}%
\fi\endgroup%
\begin{picture}(5829,744)(34,17)
\end{picture}%
\end{center}

\noindent If the announced conjecture is true, then this will prove that the first and the third representatives are identified by $\equi{0}$ and hence by $\equiv$. 
\end{ex}

\noindent The $2$-polygraphic translation of Petri nets we have built in this section has the advantage of having graphical representations that are easy to draw and interpret. However, as we have seen, the explicit way in which it handles the intrinsic commutativity of the net raises many issues we have only started to study here. The non distinction of tokens might be even worse since relations $\tau_{x,x}\equiv\id_{x\tens x}$ will create many nasty critical pairs when added to a rewriting system. However, future work will be devoted to a thorough study of these polygraphs.

The next section is devoted to a much more natural translation of Petri nets that unveils their intrinsic $3$-dimensional nature.

\section{Petri nets as $\mathbf{3}$-dimensional objects}\label{sec:3d}

\noindent In this section, we prove that Petri nets are exactly $3$-dimensional polygraphs with one cell of dimension~$0$ and no cell of dimension~$1$. The $2$-cells are the places of the net, while the $3$-cells are its transitions: there is no need of extra explicit permutation cells. This is due to a topological properties of this class of polygraphs which comes from the folkloric result of algebra, attributed to Hilton:

\begin{lem}\label{lem:algebre}
Let $M$ be a set equipped with two monoid structures $(\bullet,e)$ and $(\star,1)$ such that, for every elements~$x$, $y$, $z$ and $t$ in $M$, the relation $(x\bullet y)\star(z\bullet t)=(x\star z)\bullet(y\star t)$ holds. Then the two monoid structures are equal and commutative, which means that $e=1$ and that $x\bullet y=x\star y=y\bullet x=y\star x$.
\end{lem}

\begin{dem}
Let us start by proving the equality $e=1$. Let us apply the hypothesis with $x=t=e$ and $y=z=1$, which gives $(e\bullet 1)\star(1\bullet e)=(e\star 1)\bullet(1\star e)$. On one hand, we have $(e\bullet 1)\star(1\bullet e)=1\star 1=1$, since $e$ is a bilateral unit for $\bullet$ and since $1$ is a left (or right) unit for $\star$. But, on the other hand, $(e\star 1)\bullet(1\star e)=e\bullet e=e$, since $1$ is a bilateral unit for $\star$ and since $e$ is a left (or right) unit for $\bullet$. Hence $e=1$. 

In order to prove that both operations $\bullet$ and $\star$ are the same, let us fix two elements $x$ and $y$ in $M$. We have the following chain of equalities, using the hypothesis together with the facts that $1$ is a bilateral unit for $\star$ and for $\bullet$:
$$
x\bullet y \:=\: (x\star 1)\bullet(1\star y) \:=\: (x\bullet 1)\star(1\bullet y) \:=\: x\star y.
$$

\noindent Finally, we prove that the operation $\star$ is commutative, using the same arguments:
$$
x\star y \:=\: (1\bullet x)\star(y\bullet 1) \:=\: (1\star y)\bullet(x\star 1) \:=\: y\bullet x \:=\: y\star x.
$$
\findem\end{dem}

\begin{rem}
The proof does not use the associativity of $\bullet$ nor $\star$. It works with a set with two binary relations such that each one admits a bilateral unit.
\end{rem}

\noindent Let us translate the lemma \ref{lem:algebre} in our setting:

\begin{cor}\label{cor:algebre}
Let $\Sigma=(\ast,\Sigma_1,\Sigma_2)$ be a $2$-polygraph with one $0$-cell. Then the two compositions $\tens$ and~$\circ$ are equal and commutative on the set $\mon{\Sigma}_2(\id_{\ast},\id_{\ast})$, which is the set of all the $2$-arrows $\id_{\ast}\fl\id_{\ast}$ of the free $2$-category $\mon{\Sigma}$.
\end{cor}

\begin{dem}
On $\mon{\Sigma}_2(\id_{\ast},\id_{\ast})$ both compositions $\tens$ and $\circ$ induce a monoid structure. We already know that both structures have the same neutral element, $\id_{\id_{\ast}}$. Furthermore, the exchange relation gives, for any four $f$, $g$, $h$ and $k$ in $\mon{\Sigma}_2(\id_{\ast},\id_{\ast})$:
$$
(f\tens g)\circ(h\tens k) \:=\: (f\circ h)\tens(g\circ k).
$$

\noindent Then, one applies lemma \ref{lem:algebre} to conclude.
\findem\end{dem}

\medskip
\begin{ntt}
Let $\Sigma=(\ast,\Sigma_1,\Sigma_2)$ be a $2$-polygraph with one $0$-cell. The $1$-arrow $\id_{\ast}$ is denoted by $0$ and, by a slight abuse, so is the $2$-arrow $\id_{\id_{\ast}}$. The common restriction of $\circ$ and $\tens$ to $\mon{\Sigma}_2(0,0)$ is denoted by $+$. 
\end{ntt}

\begin{rem}
A $2$-arrow with source and target equal to $0$ is represented as a circuit with no input wire and no output wire. The proof that both compositions are equal and commutative on this kind of $2$-arrows corresponds to the following moves:
\begin{center}
\begin{picture}(0,0)%
\includegraphics{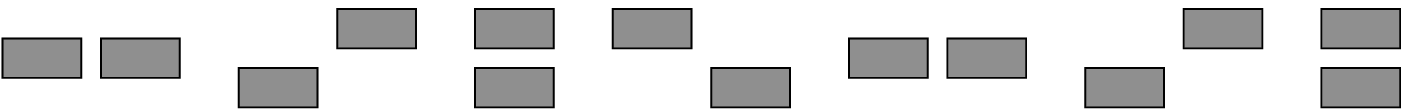}%
\end{picture}%
\setlength{\unitlength}{4144sp}%
\begingroup\makeatletter\ifx\SetFigFont\undefined%
\gdef\SetFigFont#1#2#3#4#5{%
  \reset@font\fontsize{#1}{#2pt}%
  \fontfamily{#3}\fontseries{#4}\fontshape{#5}%
  \selectfont}%
\fi\endgroup%
\begin{picture}(6414,474)(169,242)
\put(6391,299){\makebox(0,0)[b]{\smash{{\SetFigFont{10}{12.0}{\rmdefault}{\mddefault}{\updefault}{\color[rgb]{0,0,0}$g$}%
}}}}
\put(1126,434){\makebox(0,0)[b]{\smash{{\SetFigFont{10}{12.0}{\rmdefault}{\mddefault}{\updefault}{\color[rgb]{0,0,0}$=$}%
}}}}
\put(2206,434){\makebox(0,0)[b]{\smash{{\SetFigFont{10}{12.0}{\rmdefault}{\mddefault}{\updefault}{\color[rgb]{0,0,0}$=$}%
}}}}
\put(2836,434){\makebox(0,0)[b]{\smash{{\SetFigFont{10}{12.0}{\rmdefault}{\mddefault}{\updefault}{\color[rgb]{0,0,0}$=$}%
}}}}
\put(3916,434){\makebox(0,0)[b]{\smash{{\SetFigFont{10}{12.0}{\rmdefault}{\mddefault}{\updefault}{\color[rgb]{0,0,0}$=$}%
}}}}
\put(4996,434){\makebox(0,0)[b]{\smash{{\SetFigFont{10}{12.0}{\rmdefault}{\mddefault}{\updefault}{\color[rgb]{0,0,0}$=$}%
}}}}
\put(6076,434){\makebox(0,0)[b]{\smash{{\SetFigFont{10}{12.0}{\rmdefault}{\mddefault}{\updefault}{\color[rgb]{0,0,0}$=$}%
}}}}
\put(361,434){\makebox(0,0)[b]{\smash{{\SetFigFont{10}{12.0}{\rmdefault}{\mddefault}{\updefault}{\color[rgb]{0,0,0}$f$}%
}}}}
\put(811,434){\makebox(0,0)[b]{\smash{{\SetFigFont{10}{12.0}{\rmdefault}{\mddefault}{\updefault}{\color[rgb]{0,0,0}$g$}%
}}}}
\put(1441,299){\makebox(0,0)[b]{\smash{{\SetFigFont{10}{12.0}{\rmdefault}{\mddefault}{\updefault}{\color[rgb]{0,0,0}$f$}%
}}}}
\put(2521,299){\makebox(0,0)[b]{\smash{{\SetFigFont{10}{12.0}{\rmdefault}{\mddefault}{\updefault}{\color[rgb]{0,0,0}$f$}%
}}}}
\put(3601,299){\makebox(0,0)[b]{\smash{{\SetFigFont{10}{12.0}{\rmdefault}{\mddefault}{\updefault}{\color[rgb]{0,0,0}$f$}%
}}}}
\put(4681,434){\makebox(0,0)[b]{\smash{{\SetFigFont{10}{12.0}{\rmdefault}{\mddefault}{\updefault}{\color[rgb]{0,0,0}$f$}%
}}}}
\put(5761,569){\makebox(0,0)[b]{\smash{{\SetFigFont{10}{12.0}{\rmdefault}{\mddefault}{\updefault}{\color[rgb]{0,0,0}$f$}%
}}}}
\put(6391,569){\makebox(0,0)[b]{\smash{{\SetFigFont{10}{12.0}{\rmdefault}{\mddefault}{\updefault}{\color[rgb]{0,0,0}$f$}%
}}}}
\put(1891,569){\makebox(0,0)[b]{\smash{{\SetFigFont{10}{12.0}{\rmdefault}{\mddefault}{\updefault}{\color[rgb]{0,0,0}$g$}%
}}}}
\put(2521,569){\makebox(0,0)[b]{\smash{{\SetFigFont{10}{12.0}{\rmdefault}{\mddefault}{\updefault}{\color[rgb]{0,0,0}$g$}%
}}}}
\put(3151,569){\makebox(0,0)[b]{\smash{{\SetFigFont{10}{12.0}{\rmdefault}{\mddefault}{\updefault}{\color[rgb]{0,0,0}$g$}%
}}}}
\put(4231,434){\makebox(0,0)[b]{\smash{{\SetFigFont{10}{12.0}{\rmdefault}{\mddefault}{\updefault}{\color[rgb]{0,0,0}$g$}%
}}}}
\put(5311,299){\makebox(0,0)[b]{\smash{{\SetFigFont{10}{12.0}{\rmdefault}{\mddefault}{\updefault}{\color[rgb]{0,0,0}$g$}%
}}}}
\end{picture}%
\end{center}

\noindent Explicitely:
$$
f\tens g \:=\: (f\tens 0)\circ(0\tens g) \:=\: f\circ g \:=\: (0\tens f)\circ(g\tens 0) \:=\: g\tens f \:=\: (g\tens 0)\circ(0\tens f) \:=\: g\circ f.
$$

\noindent This means that such a special $2$-arrow can turn around another one: there is no wire, hence no limitation to their homeomorphic movement.
\end{rem}

\noindent Using corollary \ref{prop:comm}, we give a polygraphic description of the free commutative monoid generated by a given set:

\begin{prop}\label{prop:comm}
Let $\Sigma$ be a $2$-polygraph of the form $(\ast,\emptyset,\Sigma_2)$. Then the set $\mon{\Sigma}_2(0,0)$ contains all the~$2$-arrows of $\mon{\Sigma}$ and, equipped with the structure $(+,0)$, is isomorphic to the free commutative monoid~$[\Sigma_2]$ generated by $\Sigma_2$.
\end{prop}

\begin{dem}
Since there is one $0$-cell and no $1$-cell in the $2$-polygraph $\Sigma$, the only $1$-arrow of the free $2$-category~$\mon{\Sigma}$ is $\id_{\ast}=0$: indeed, there is only one path in the graph $(\ast,\emptyset)$ with one object and no arrow, the empty one. Hence, every $2$-arrow of $\mon{\Sigma}$ starts and ends at $0$. 

By application of corollary \ref{cor:algebre}, we know that both compositions $\circ$ and $\tens$ are equal and commutative, so that $(\mon{\Sigma}_2,+,0)$ is a commutative monoid. Furthermore, each element of $\Sigma_2$ is represented in $\mon{\Sigma}_2$: this inclusion induces a unique monoid morphism from $[\Sigma_2]$ to $\mon{\Sigma}_2$. This morphism is surjective, since every $2$-arrow of $\mon{\Sigma}$ is built from $2$-cells (elements of $\Sigma_2$) using only the operations $\tens$ and $\circ$, both equal to $+$. Hence, every $2$-arrow $f$ of $\mon{\Sigma}$ admits a decomposition:
$$
f \:=\: \sum_{x\in\Sigma_2} f(x).x,
$$

\noindent where the $f(x)$ are natural numbers. In order to conclude the proof, one must prove that this decomposition is unique. Let us assume that $f$ has another decomposition:
$$
f \:=\: \sum_{x\in\Sigma_2} f'(x).x.
$$

\noindent Let us fix a $2$-cell $x\in\Sigma_2$ and assume that $f(x)=f'(x)+k$, with $k$ a natural number. Then:
$$
f-f'(x).x \:=\: \sum_{y\neq x} f(y).y + k.x \:=\: \sum_{y\neq x} f'(y).y.
$$

\noindent Hence, in the first decomposition of $f-f'(x).x$, there are $k$ copies of the $2$-cell $x$, but there are no in the other. However, in a free $2$-category, two arrows are equal if and only if they differ only by a limited number of applications of the rules of associativity, units and exchange for $\circ$ and $\tens$: all these operations leave the number of generating $2$-cells $x$ unchanged. Hence $k=0$ and $f(x)=f'(x)$. Finally, there are only a finite number of $x$ such that $f(x)\neq 0$: an induction on this number conludes the proof.

\findem\end{dem}

\noindent Now, we have a correspondance between the elements of $[X]$ and the $2$-arrows of the free $2$-category generated by $(\ast,\emptyset,X)$. Then, transitions of a Petri net, through their rewriting representation, are translated as $\mathit{3}$-\emph{cells} in a $\mathit{3}$-\emph{polygraph}.

\begin{defn}
A \emph{$\mathit{3}$-polygraph} is a family $\Sigma=(\Sigma_0,\Sigma_1,\Sigma_2,\Sigma_3)$ of sets, equipped with an additional structure of $2$-polygraph on $(\Sigma_0,\Sigma_1,\Sigma_2)$ and with a graph structure $s_2,t_2:\Sigma_3\fl\mon{\Sigma}_2$ such that:
$$
s_1\circ s_2 \:=\: s_1\circ t_2 \et t_1\circ s_2 \:=\: t_1\circ t_2.
$$
\end{defn}

\begin{rem}
Usually, the $3$-cells are seen as \emph{directed volumes} between parallel circuits (circuits with the same $1$-source and the same $1$-target).
\end{rem}

\noindent Let us formalize the translation from Petri nets into $3$-polygraphs:

\begin{defn}
Let $(X,R)$ be a Petri net. The \emph{$\mathit{3}$-polygraph associated to $(X,R)$}, denoted by $\Sigma^3(X,R)$, is the $3$-polygraph $(\ast,\emptyset,X,R)$, where each rewriting rule $\alpha=(a,b)$ is seen as a $3$-cell with $2$-source the circuit representing $a$ and $2$-target the circuit representing $b$.

Conversely, let $\Sigma=(\ast,\emptyset,\Sigma_2,\Sigma_3)$ be a $3$-polygraph with one $0$-cell and no $1$-cell. Its associated Petri net is the pair $\Nr(\Sigma)=(\Sigma_2,\Sigma_3)$.
\end{defn}

\noindent In order to compare a Petri net and its associated $3$-polygraph, a notion of reduction graph is defined, which conveys the idea of reduction under a context - see [Guiraud 2004(T)] for a study of contexts for circuits:

\begin{defn}
Let $\Sigma=(\ast,\Sigma_1,\Sigma_2,\Sigma_3)$ be a $3$-polygraph with one $0$-cell. Its \emph{associated reduction graph} is the graph $G(\Sigma)$ defined this way:
\begin{enumerate}
\item[0.] The objects of $G(\Sigma)$ are the $2$-arrows of $\mon{\Sigma}_2$.
\item[1.] The arrows of $G(\Sigma)$ from $u$ to $v$ are all the triples $(f,\alpha,g)$, made of two $2$-arrows $f$ and $g$ of $\mon{\Sigma}_2$ and one $3$-cell $\alpha$ of $\Sigma_3$, such that the two following equalities are defined and hold:
\begin{center}
\begin{picture}(0,0)%
\includegraphics{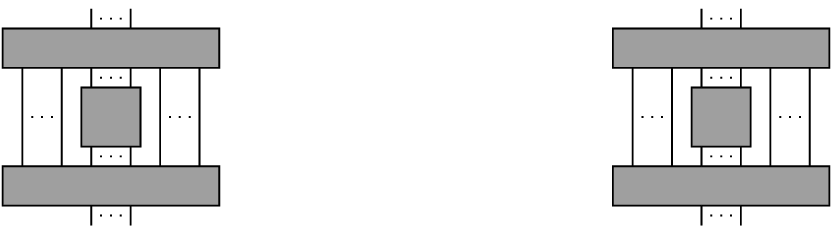}%
\end{picture}%
\setlength{\unitlength}{4144sp}%
\begingroup\makeatletter\ifx\SetFigFont\undefined%
\gdef\SetFigFont#1#2#3#4#5{%
  \reset@font\fontsize{#1}{#2pt}%
  \fontfamily{#3}\fontseries{#4}\fontshape{#5}%
  \selectfont}%
\fi\endgroup%
\begin{picture}(4447,1014)(79,-253)
\put(4411,209){\makebox(0,0)[b]{\smash{{\SetFigFont{10}{12.0}{\rmdefault}{\mddefault}{\updefault}{\color[rgb]{0,0,0}$v$}%
}}}}
\put(586,524){\makebox(0,0)[b]{\smash{{\SetFigFont{10}{12.0}{\rmdefault}{\mddefault}{\updefault}{\color[rgb]{0,0,0}$f$}%
}}}}
\put(3376,524){\makebox(0,0)[b]{\smash{{\SetFigFont{10}{12.0}{\rmdefault}{\mddefault}{\updefault}{\color[rgb]{0,0,0}$f$}%
}}}}
\put(586,-106){\makebox(0,0)[b]{\smash{{\SetFigFont{10}{12.0}{\rmdefault}{\mddefault}{\updefault}{\color[rgb]{0,0,0}$g$}%
}}}}
\put(3376,-106){\makebox(0,0)[b]{\smash{{\SetFigFont{10}{12.0}{\rmdefault}{\mddefault}{\updefault}{\color[rgb]{0,0,0}$g$}%
}}}}
\put(586,209){\makebox(0,0)[b]{\smash{{\SetFigFont{10}{12.0}{\rmdefault}{\mddefault}{\updefault}{\color[rgb]{0,0,0}$s_2\alpha$}%
}}}}
\put(3376,209){\makebox(0,0)[b]{\smash{{\SetFigFont{10}{12.0}{\rmdefault}{\mddefault}{\updefault}{\color[rgb]{0,0,0}$t_2\alpha$}%
}}}}
\put(1351,209){\makebox(0,0)[b]{\smash{{\SetFigFont{10}{12.0}{\rmdefault}{\mddefault}{\updefault}{\color[rgb]{0,0,0}$=$}%
}}}}
\put(1621,209){\makebox(0,0)[b]{\smash{{\SetFigFont{10}{12.0}{\rmdefault}{\mddefault}{\updefault}{\color[rgb]{0,0,0}$u$}%
}}}}
\put(4141,209){\makebox(0,0)[b]{\smash{{\SetFigFont{10}{12.0}{\rmdefault}{\mddefault}{\updefault}{\color[rgb]{0,0,0}$=$}%
}}}}
\end{picture}%
\end{center}
\noindent These triples are considered \emph{modulo} the following \emph{deformation} equations, given for every possible $2$-arrows~$f$,~$g$ and~$h$ and $3$-cell $\alpha$:
\begin{center}
\begin{picture}(0,0)%
\includegraphics{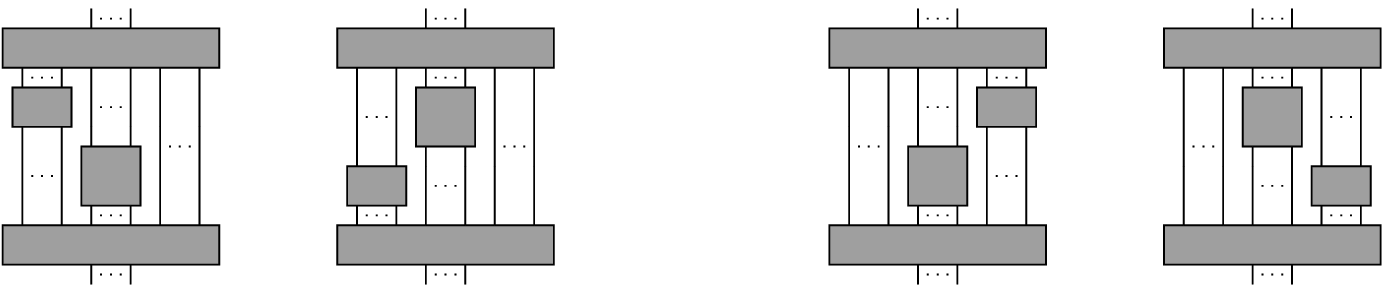}%
\end{picture}%
\setlength{\unitlength}{4144sp}%
\begingroup\makeatletter\ifx\SetFigFont\undefined%
\gdef\SetFigFont#1#2#3#4#5{%
  \reset@font\fontsize{#1}{#2pt}%
  \fontfamily{#3}\fontseries{#4}\fontshape{#5}%
  \selectfont}%
\fi\endgroup%
\begin{picture}(6324,1284)(79,-253)
\put(5896,479){\makebox(0,0)[b]{\smash{{\SetFigFont{10}{12.0}{\rmdefault}{\mddefault}{\updefault}{\color[rgb]{0,0,0}$\alpha$}%
}}}}
\put(1351,389){\makebox(0,0)[b]{\smash{{\SetFigFont{10}{12.0}{\rmdefault}{\mddefault}{\updefault}{\color[rgb]{0,0,0}$=$}%
}}}}
\put(5131,389){\makebox(0,0)[b]{\smash{{\SetFigFont{10}{12.0}{\rmdefault}{\mddefault}{\updefault}{\color[rgb]{0,0,0}$=$}%
}}}}
\put(586,794){\makebox(0,0)[b]{\smash{{\SetFigFont{10}{12.0}{\rmdefault}{\mddefault}{\updefault}{\color[rgb]{0,0,0}$f$}%
}}}}
\put(2116,794){\makebox(0,0)[b]{\smash{{\SetFigFont{10}{12.0}{\rmdefault}{\mddefault}{\updefault}{\color[rgb]{0,0,0}$f$}%
}}}}
\put(4366,794){\makebox(0,0)[b]{\smash{{\SetFigFont{10}{12.0}{\rmdefault}{\mddefault}{\updefault}{\color[rgb]{0,0,0}$f$}%
}}}}
\put(5896,794){\makebox(0,0)[b]{\smash{{\SetFigFont{10}{12.0}{\rmdefault}{\mddefault}{\updefault}{\color[rgb]{0,0,0}$f$}%
}}}}
\put(586,-106){\makebox(0,0)[b]{\smash{{\SetFigFont{10}{12.0}{\rmdefault}{\mddefault}{\updefault}{\color[rgb]{0,0,0}$g$}%
}}}}
\put(2116,-106){\makebox(0,0)[b]{\smash{{\SetFigFont{10}{12.0}{\rmdefault}{\mddefault}{\updefault}{\color[rgb]{0,0,0}$g$}%
}}}}
\put(4681,524){\makebox(0,0)[b]{\smash{{\SetFigFont{10}{12.0}{\rmdefault}{\mddefault}{\updefault}{\color[rgb]{0,0,0}$h$}%
}}}}
\put(6211,164){\makebox(0,0)[b]{\smash{{\SetFigFont{10}{12.0}{\rmdefault}{\mddefault}{\updefault}{\color[rgb]{0,0,0}$h$}%
}}}}
\put(271,524){\makebox(0,0)[b]{\smash{{\SetFigFont{10}{12.0}{\rmdefault}{\mddefault}{\updefault}{\color[rgb]{0,0,0}$h$}%
}}}}
\put(1801,164){\makebox(0,0)[b]{\smash{{\SetFigFont{10}{12.0}{\rmdefault}{\mddefault}{\updefault}{\color[rgb]{0,0,0}$h$}%
}}}}
\put(4366,-106){\makebox(0,0)[b]{\smash{{\SetFigFont{10}{12.0}{\rmdefault}{\mddefault}{\updefault}{\color[rgb]{0,0,0}$g$}%
}}}}
\put(5896,-106){\makebox(0,0)[b]{\smash{{\SetFigFont{10}{12.0}{\rmdefault}{\mddefault}{\updefault}{\color[rgb]{0,0,0}$g$}%
}}}}
\put(586,209){\makebox(0,0)[b]{\smash{{\SetFigFont{10}{12.0}{\rmdefault}{\mddefault}{\updefault}{\color[rgb]{0,0,0}$\alpha$}%
}}}}
\put(2116,479){\makebox(0,0)[b]{\smash{{\SetFigFont{10}{12.0}{\rmdefault}{\mddefault}{\updefault}{\color[rgb]{0,0,0}$\alpha$}%
}}}}
\put(4366,209){\makebox(0,0)[b]{\smash{{\SetFigFont{10}{12.0}{\rmdefault}{\mddefault}{\updefault}{\color[rgb]{0,0,0}$\alpha$}%
}}}}
\end{picture}%
\end{center}
\end{enumerate}

\noindent A triple $(f,\alpha,g)$ is denoted by $g\circ\alpha\circ f$, with $(\circ f)$ and/or $(g\circ)$ dropped when $f$ and/or $g$ is an identity. We denote by $\star$ the composition of the free category $\mon{G(\Sigma)}$, with $A\star B$ standing $A$ followed by $B$.
\end{defn}

\noindent Once again, the reduction graph is not the natural object one associates to a $3$-polygraph: we prefer the $3$-category it generates. We give a formal definition and, then, its underlying graphical intuition.

\begin{defn}
Let $\Sigma=(\ast,\Sigma_1,\Sigma_2,\Sigma_3)$ be a $3$-polygraph with one $0$-cell. The \emph{free $\mathit{3}$-category generated by $\Sigma$} is denoted by $\mon{\Sigma}$ and is made of the $0$, $1$ and $2$-arrows of $\Sigma$, together with a family of \emph{$\mathit{3}$-arrows} which are the paths of the reduction graph $G(\Sigma)$ \emph{modulo} the congruence $\equi{\Sigma}$ generated by the following exchange relations:
$$
\begin{array}{c c c}
(A\tens s_2(B))\star(t_2(A)\tens B) 
&\quad\equi{02}\quad&
(s_2(A)\tens B)\star(A\tens t_2(B)), \\
(B\circ s_2(A))\star(t_2(B)\circ A)
&\quad\equi{12}\quad&
(s_2(B)\circ A)\star(B\circ t_2(A)).
\end{array}
$$

\noindent These equations allow one to extend the two compositions $\tens$ and $\circ$ on equivalence classes of paths in the graph $G(\Sigma)$, with $A\tens B$ being given by either side of the relation $\equi{02}$ and $B\circ A$ by either side of $\equi{12}$.
\end{defn}

\begin{rem}
Let us give a more graphical account of the free $3$-category $\mon{\Sigma}$ generated by a $3$-polygraph $\Sigma$. Its $3$-arrows are generated by the $3$-cells of $\Sigma$ seen as blocks:
\begin{center}
\begin{picture}(0,0)%
\includegraphics{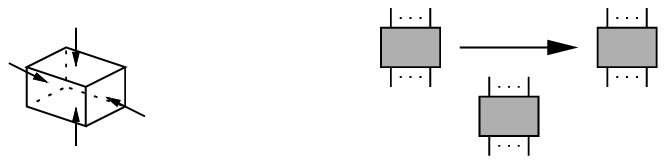}%
\end{picture}%
\setlength{\unitlength}{4144sp}%
\begingroup\makeatletter\ifx\SetFigFont\undefined%
\gdef\SetFigFont#1#2#3#4#5{%
  \reset@font\fontsize{#1}{#2pt}%
  \fontfamily{#3}\fontseries{#4}\fontshape{#5}%
  \selectfont}%
\fi\endgroup%
\begin{picture}(3225,868)(118,-146)
\put(1126, 29){\makebox(0,0)[b]{\smash{{\SetFigFont{10}{12.0}{\rmdefault}{\mddefault}{\updefault}{\color[rgb]{0,0,0}$g$}%
}}}}
\put(676,614){\makebox(0,0)[b]{\smash{{\SetFigFont{10}{12.0}{\rmdefault}{\mddefault}{\updefault}{\color[rgb]{0,0,0}$m$}%
}}}}
\put(676,-106){\makebox(0,0)[b]{\smash{{\SetFigFont{10}{12.0}{\rmdefault}{\mddefault}{\updefault}{\color[rgb]{0,0,0}$n$}%
}}}}
\put(2206,434){\makebox(0,0)[b]{\smash{{\SetFigFont{10}{12.0}{\rmdefault}{\mddefault}{\updefault}$f$}}}}
\put(3196,434){\makebox(0,0)[b]{\smash{{\SetFigFont{10}{12.0}{\rmdefault}{\mddefault}{\updefault}$g$}}}}
\put(2656,119){\makebox(0,0)[b]{\smash{{\SetFigFont{10}{12.0}{\rmdefault}{\mddefault}{\updefault}$A$}}}}
\put(226,389){\makebox(0,0)[b]{\smash{{\SetFigFont{10}{12.0}{\rmdefault}{\mddefault}{\updefault}{\color[rgb]{0,0,0}$f$}%
}}}}
\end{picture}%
\end{center}

\noindent On these generators, one can use the three following constructors, called compositions:
\begin{center}
\begin{picture}(0,0)%
\includegraphics{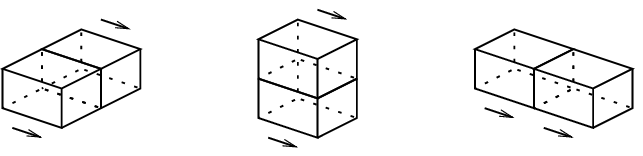}%
\end{picture}%
\setlength{\unitlength}{4144sp}%
\begingroup\makeatletter\ifx\SetFigFont\undefined%
\gdef\SetFigFont#1#2#3#4#5{%
  \reset@font\fontsize{#1}{#2pt}%
  \fontfamily{#3}\fontseries{#4}\fontshape{#5}%
  \selectfont}%
\fi\endgroup%
\begin{picture}(2904,1249)(169,-470)
\put(2701,-421){\makebox(0,0)[b]{\smash{{\SetFigFont{10}{12.0}{\rmdefault}{\mddefault}{\updefault}{\color[rgb]{0,0,0}$A\star B$}%
}}}}
\put(766,614){\makebox(0,0)[b]{\smash{{\SetFigFont{10}{12.0}{\rmdefault}{\mddefault}{\updefault}{\color[rgb]{0,0,0}$B$}%
}}}}
\put(226,-106){\makebox(0,0)[b]{\smash{{\SetFigFont{10}{12.0}{\rmdefault}{\mddefault}{\updefault}{\color[rgb]{0,0,0}$A$}%
}}}}
\put(1396,-151){\makebox(0,0)[b]{\smash{{\SetFigFont{10}{12.0}{\rmdefault}{\mddefault}{\updefault}{\color[rgb]{0,0,0}$B$}%
}}}}
\put(2386,-16){\makebox(0,0)[b]{\smash{{\SetFigFont{10}{12.0}{\rmdefault}{\mddefault}{\updefault}{\color[rgb]{0,0,0}$A$}%
}}}}
\put(2656,-106){\makebox(0,0)[b]{\smash{{\SetFigFont{10}{12.0}{\rmdefault}{\mddefault}{\updefault}{\color[rgb]{0,0,0}$B$}%
}}}}
\put(1756,659){\makebox(0,0)[b]{\smash{{\SetFigFont{10}{12.0}{\rmdefault}{\mddefault}{\updefault}{\color[rgb]{0,0,0}$A$}%
}}}}
\put(496,-421){\makebox(0,0)[b]{\smash{{\SetFigFont{10}{12.0}{\rmdefault}{\mddefault}{\updefault}{\color[rgb]{0,0,0}$A\tens B$}%
}}}}
\put(1576,-421){\makebox(0,0)[b]{\smash{{\SetFigFont{10}{12.0}{\rmdefault}{\mddefault}{\updefault}{\color[rgb]{0,0,0}$B\circ A$}%
}}}}
\end{picture}%
\end{center}

\noindent If they are sliced, these compositions appear this way:
\begin{center}
\begin{picture}(0,0)%
\includegraphics{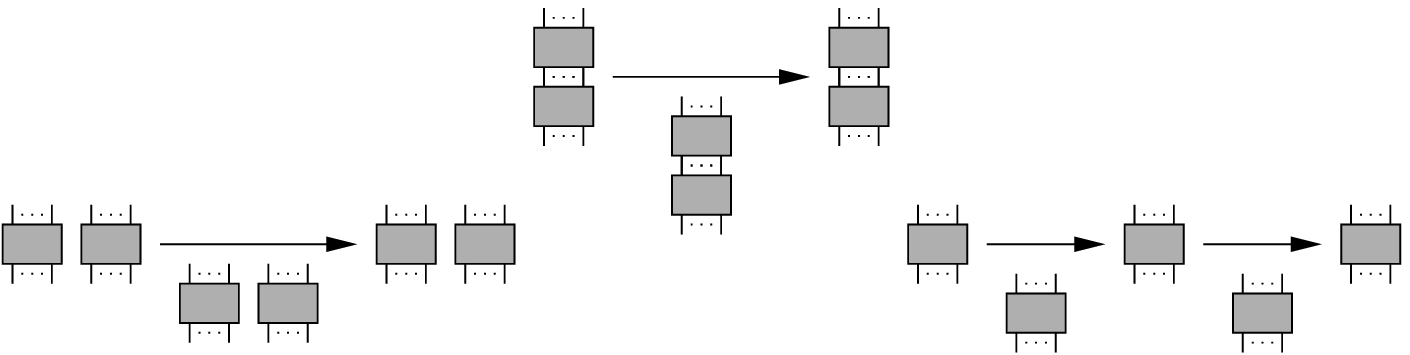}%
\end{picture}%
\setlength{\unitlength}{4144sp}%
\begingroup\makeatletter\ifx\SetFigFont\undefined%
\gdef\SetFigFont#1#2#3#4#5{%
  \reset@font\fontsize{#1}{#2pt}%
  \fontfamily{#3}\fontseries{#4}\fontshape{#5}%
  \selectfont}%
\fi\endgroup%
\begin{picture}(6414,1599)(79,-838)
\put(2656,254){\makebox(0,0)[b]{\smash{{\SetFigFont{10}{12.0}{\rmdefault}{\mddefault}{\updefault}$g$}}}}
\put(4006,524){\makebox(0,0)[b]{\smash{{\SetFigFont{10}{12.0}{\rmdefault}{\mddefault}{\updefault}$f'$}}}}
\put(4006,254){\makebox(0,0)[b]{\smash{{\SetFigFont{10}{12.0}{\rmdefault}{\mddefault}{\updefault}$g'$}}}}
\put(4366,-376){\makebox(0,0)[b]{\smash{{\SetFigFont{10}{12.0}{\rmdefault}{\mddefault}{\updefault}$f$}}}}
\put(5356,-376){\makebox(0,0)[b]{\smash{{\SetFigFont{10}{12.0}{\rmdefault}{\mddefault}{\updefault}$g$}}}}
\put(6346,-376){\makebox(0,0)[b]{\smash{{\SetFigFont{10}{12.0}{\rmdefault}{\mddefault}{\updefault}$h$}}}}
\put(1036,-646){\makebox(0,0)[b]{\smash{{\SetFigFont{10}{12.0}{\rmdefault}{\mddefault}{\updefault}$A$}}}}
\put(1396,-646){\makebox(0,0)[b]{\smash{{\SetFigFont{10}{12.0}{\rmdefault}{\mddefault}{\updefault}$B$}}}}
\put(3286,119){\makebox(0,0)[b]{\smash{{\SetFigFont{10}{12.0}{\rmdefault}{\mddefault}{\updefault}$A$}}}}
\put(3286,-151){\makebox(0,0)[b]{\smash{{\SetFigFont{10}{12.0}{\rmdefault}{\mddefault}{\updefault}$B$}}}}
\put(4816,-691){\makebox(0,0)[b]{\smash{{\SetFigFont{10}{12.0}{\rmdefault}{\mddefault}{\updefault}$A$}}}}
\put(5851,-691){\makebox(0,0)[b]{\smash{{\SetFigFont{10}{12.0}{\rmdefault}{\mddefault}{\updefault}$B$}}}}
\put(226,-376){\makebox(0,0)[b]{\smash{{\SetFigFont{10}{12.0}{\rmdefault}{\mddefault}{\updefault}$f$}}}}
\put(586,-376){\makebox(0,0)[b]{\smash{{\SetFigFont{10}{12.0}{\rmdefault}{\mddefault}{\updefault}$g$}}}}
\put(1936,-376){\makebox(0,0)[b]{\smash{{\SetFigFont{10}{12.0}{\rmdefault}{\mddefault}{\updefault}$f'$}}}}
\put(2296,-376){\makebox(0,0)[b]{\smash{{\SetFigFont{10}{12.0}{\rmdefault}{\mddefault}{\updefault}$g'$}}}}
\put(2656,524){\makebox(0,0)[b]{\smash{{\SetFigFont{10}{12.0}{\rmdefault}{\mddefault}{\updefault}$f$}}}}
\end{picture}%
\end{center}

\vfill\pagebreak
\noindent All the constructions are identified modulo the following moves:
\begin{center}
\begin{picture}(0,0)%
\includegraphics{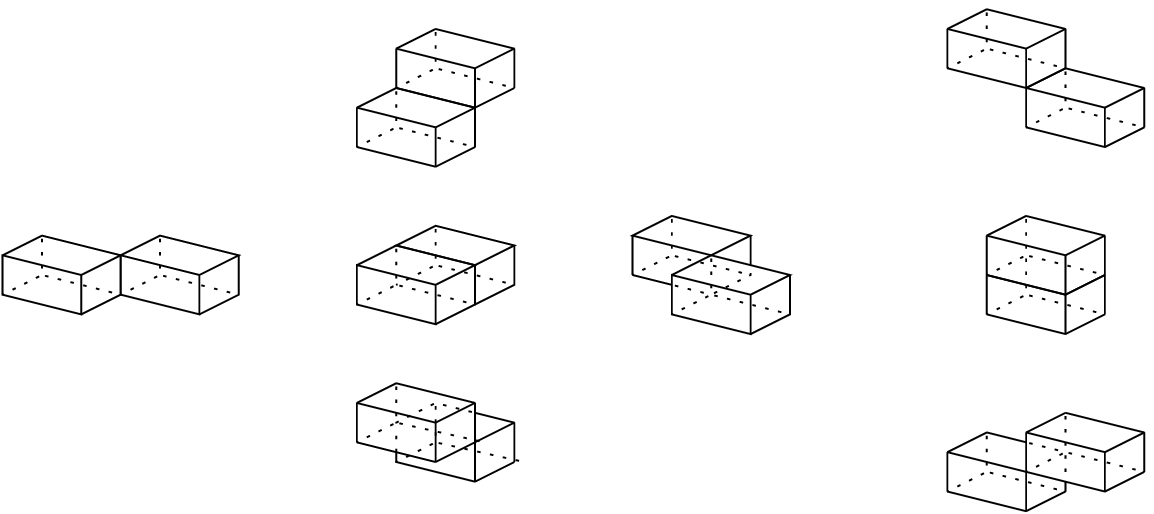}%
\end{picture}%
\setlength{\unitlength}{4144sp}%
\begingroup\makeatletter\ifx\SetFigFont\undefined%
\gdef\SetFigFont#1#2#3#4#5{%
  \reset@font\fontsize{#1}{#2pt}%
  \fontfamily{#3}\fontseries{#4}\fontshape{#5}%
  \selectfont}%
\fi\endgroup%
\begin{picture}(5244,2319)(169,-1513)
\put(2791,-466){\makebox(0,0)[b]{\smash{{\SetFigFont{10}{12.0}{\rmdefault}{\mddefault}{\updefault}{\color[rgb]{0,0,0}$\equi{02}$}%
}}}}
\put(2161,-106){\makebox(0,0)[b]{\smash{{\SetFigFont{10}{12.0}{\rmdefault}{\mddefault}{\updefault}{\color[rgb]{0,0,0}$\equi{01}$}%
}}}}
\put(2161,-826){\makebox(0,0)[b]{\smash{{\SetFigFont{10}{12.0}{\rmdefault}{\mddefault}{\updefault}{\color[rgb]{0,0,0}$\equi{01}$}%
}}}}
\put(4951,-61){\makebox(0,0)[b]{\smash{{\SetFigFont{10}{12.0}{\rmdefault}{\mddefault}{\updefault}{\color[rgb]{0,0,0}$\equi{12}$}%
}}}}
\put(4951,-916){\makebox(0,0)[b]{\smash{{\SetFigFont{10}{12.0}{\rmdefault}{\mddefault}{\updefault}{\color[rgb]{0,0,0}$\equi{12}$}%
}}}}
\put(1531,-466){\makebox(0,0)[b]{\smash{{\SetFigFont{10}{12.0}{\rmdefault}{\mddefault}{\updefault}{\color[rgb]{0,0,0}$\equi{02}$}%
}}}}
\end{picture}%
\end{center}

\noindent This picture contains three families of moves, one for each exchange relation $\equi{02}$, $\equi{12}$ and $\equi{01}$, where the relation $\equi{01}$ is induced by the deformation relations and the other two exchange relations.
\end{rem}

\begin{rem}
In the case of a $3$-polygraph $\Sigma$ with one $0$-cell and no $1$-cell, there are only two ways to compose $3$-arrows, namely $+$ and $\star$, since $\circ$ and $\tens$ are the same and denoted by $+$. As a consequence, there is only one family of exchange relations:
$$
(A+s_2(B))\star(t_2(A)+B) \quad\equi{\Sigma}\quad (s_2(A)+B)\star(A+t_2(B)). 
$$

\end{rem}

\noindent We prove that the reduction graphs of a Petri net and of its associated $3$-polygraph are the same. Moreover, the $3$-arrows of the $3$-category generated by the latter are exactly the equivalence classes of Petri net reductions \emph{modulo} the congruence relation we have defined on them.

\begin{thm}\label{th:3d}
Let $(X,R)$ be a commutative word rewriting system. Then $\Nr(\Sigma^3(X,R))=(X,R)$ and the graphs $G(\Sigma^3(X,R))$ and $G(X,R)$ are isomorphic. Furthermore, this isomorphism identifies the congruences $\equi{(X,R)}$ and $\equi{\Sigma^3(X,R)}$. Conversely, given any $3$-polygraph $\Sigma=(\ast,\emptyset,\Sigma_2,\Sigma_3)$, the equality $\Sigma^3(\Nr(\Sigma))=\Sigma$ holds and the graphs $G(\Sigma)$ and $G(\Sigma^3(\Nr(\Sigma)))$ are isomorphic.  Furthermore, this isomorphism identifies the congruences $\equi{\Sigma}$ and $\equi{\Nr(\Sigma)}$.
\end{thm}

\begin{dem}
Let us fix a Petri net $(X,R)$. The equality $\Nr(\Sigma^3(X,R))=(X,R)$ is immediate. The objects of both graphs $G(X,R)$ and of $G(\Sigma^3(X,R))$ are the same: the elements of the free commutative monoid~$[X]$. 

Then, the arrows from $u$ to $v$ in $G(X,R)$ are the $c+\alpha$, made of an element $c$ of $[X]$ and a rule $\alpha$ in $R$, such that $u=c+s(\alpha)$ and $v=c+t(\alpha)$. To such an arrow $c+\alpha$, we associate the arrow $\phi(c+\alpha)=(c,\alpha,0)$ in $G(\Sigma^3(X,R))$. 

Conversely, let us consider an arrow $(f,\alpha,g)$ in $G(\Sigma^3(X,R))$. Let us prove graphically that $(f,\alpha,g)=(f+g,\alpha,0)$, using the fact that all the $2$-arrows of $\mon{\Sigma^3(X,R)}_2$ have source and target $0$:
\begin{center}
\begin{picture}(0,0)%
\includegraphics{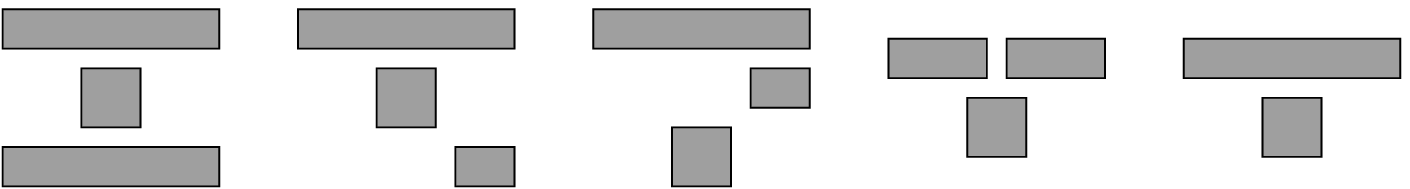}%
\end{picture}%
\setlength{\unitlength}{4144sp}%
\begingroup\makeatletter\ifx\SetFigFont\undefined%
\gdef\SetFigFont#1#2#3#4#5{%
  \reset@font\fontsize{#1}{#2pt}%
  \fontfamily{#3}\fontseries{#4}\fontshape{#5}%
  \selectfont}%
\fi\endgroup%
\begin{picture}(6414,834)(79,-163)
\put(4366,389){\makebox(0,0)[b]{\smash{{\SetFigFont{10}{12.0}{\rmdefault}{\mddefault}{\updefault}{\color[rgb]{0,0,0}$f$}%
}}}}
\put(586,209){\makebox(0,0)[b]{\smash{{\SetFigFont{10}{12.0}{\rmdefault}{\mddefault}{\updefault}{\color[rgb]{0,0,0}$\alpha$}%
}}}}
\put(586,524){\makebox(0,0)[b]{\smash{{\SetFigFont{10}{12.0}{\rmdefault}{\mddefault}{\updefault}{\color[rgb]{0,0,0}$f$}%
}}}}
\put(586,-106){\makebox(0,0)[b]{\smash{{\SetFigFont{10}{12.0}{\rmdefault}{\mddefault}{\updefault}{\color[rgb]{0,0,0}$g$}%
}}}}
\put(1936,209){\makebox(0,0)[b]{\smash{{\SetFigFont{10}{12.0}{\rmdefault}{\mddefault}{\updefault}{\color[rgb]{0,0,0}$\alpha$}%
}}}}
\put(3286,-61){\makebox(0,0)[b]{\smash{{\SetFigFont{10}{12.0}{\rmdefault}{\mddefault}{\updefault}{\color[rgb]{0,0,0}$\alpha$}%
}}}}
\put(4636, 74){\makebox(0,0)[b]{\smash{{\SetFigFont{10}{12.0}{\rmdefault}{\mddefault}{\updefault}{\color[rgb]{0,0,0}$\alpha$}%
}}}}
\put(5986,389){\makebox(0,0)[b]{\smash{{\SetFigFont{10}{12.0}{\rmdefault}{\mddefault}{\updefault}{\color[rgb]{0,0,0}$f+g$}%
}}}}
\put(5986, 74){\makebox(0,0)[b]{\smash{{\SetFigFont{10}{12.0}{\rmdefault}{\mddefault}{\updefault}{\color[rgb]{0,0,0}$\alpha$}%
}}}}
\put(1261,209){\makebox(0,0)[b]{\smash{{\SetFigFont{10}{12.0}{\rmdefault}{\mddefault}{\updefault}{\color[rgb]{0,0,0}$=$}%
}}}}
\put(2611,209){\makebox(0,0)[b]{\smash{{\SetFigFont{10}{12.0}{\rmdefault}{\mddefault}{\updefault}{\color[rgb]{0,0,0}$=$}%
}}}}
\put(3961,209){\makebox(0,0)[b]{\smash{{\SetFigFont{10}{12.0}{\rmdefault}{\mddefault}{\updefault}{\color[rgb]{0,0,0}$=$}%
}}}}
\put(5311,209){\makebox(0,0)[b]{\smash{{\SetFigFont{10}{12.0}{\rmdefault}{\mddefault}{\updefault}{\color[rgb]{0,0,0}$=$}%
}}}}
\put(1936,524){\makebox(0,0)[b]{\smash{{\SetFigFont{10}{12.0}{\rmdefault}{\mddefault}{\updefault}{\color[rgb]{0,0,0}$f$}%
}}}}
\put(2296,-106){\makebox(0,0)[b]{\smash{{\SetFigFont{10}{12.0}{\rmdefault}{\mddefault}{\updefault}{\color[rgb]{0,0,0}$g$}%
}}}}
\put(3286,524){\makebox(0,0)[b]{\smash{{\SetFigFont{10}{12.0}{\rmdefault}{\mddefault}{\updefault}{\color[rgb]{0,0,0}$f$}%
}}}}
\put(3646,254){\makebox(0,0)[b]{\smash{{\SetFigFont{10}{12.0}{\rmdefault}{\mddefault}{\updefault}{\color[rgb]{0,0,0}$g$}%
}}}}
\put(4906,389){\makebox(0,0)[b]{\smash{{\SetFigFont{10}{12.0}{\rmdefault}{\mddefault}{\updefault}{\color[rgb]{0,0,0}$g$}%
}}}}
\end{picture}%
\end{center}

\noindent Let us denote by $\psi$ the map that sends each $(f,\alpha,g)$ onto $f+g+\alpha$ and let us check that $\psi$ is an inverse for $\phi$:
$$
\psi\circ\phi(c+\alpha) \:=\: \psi(c,\alpha,0) \:=\: c+0+\alpha \:=\: c+\alpha.
$$

\noindent And:
$$
\phi\circ\psi(f,\alpha,g) \:=\: \phi(f+g+\alpha) \:=\: (f+g,\alpha,0) \:=\: (f,\alpha,g).
$$

\noindent Let us prove that $\phi(\equi{(X,R)})$ is included into $\equi{\Sigma^3(X,R)}$. For that, we fix $c$ in $[X]$ and $\alpha$, $\beta$ in $\Sigma_3$. Then $\phi$ sends the following square of $G(X,R)$ 
$$
\xymatrix
{
	c + s(\alpha) + s(\beta) \ar[rr]^-{(c+s(\beta))+\alpha} \ar[dd]_-{(c+s(\alpha))+\beta}
	&&
	c + t(\alpha) + s(\beta) \ar[dd]^-{(c+t(\alpha))+\beta} 
	\\ \\
	c + s(\alpha) + t(\beta) \ar[rr]_-{(c+t(\beta))+\alpha}
	&&
	c + t(\alpha) + t(\beta)
}
$$

\noindent onto the following square of $G(\Sigma^3(X,R))$:
$$
\xymatrix
{
	c + s(\alpha) + s(\beta) \ar[rr]^-{(c+s(\beta),\alpha,0)} \ar[dd]_-{(c+s(\alpha),\beta,0)}
	&&
	c + t(\alpha) + s(\beta) \ar[dd]^-{(c+t(\alpha),\beta,0)} 
	\\ \\
	c + s(\alpha) + t(\beta) \ar[rr]_-{(c+t(\beta),\alpha,0)}
	&&
	c + t(\alpha) + t(\beta).
}
$$

\noindent Using the already-known properties of $G(\Sigma^3(X,R))$, we get the following two equations:
$$
\left\{
\begin{array}{c c c}
(c+s(\beta),\alpha,0)\star(c+t(\alpha),\beta,0) &\:=\:&
c+\big((\alpha+s(\beta))\star(t(\alpha)+\beta)\big), \\
(c+s(\alpha),\beta,0)\star(c+t(\beta),\alpha,0) &\:=\:&
c+\big((s(\alpha)+\beta)\star(\alpha+t(\beta))\big).
\end{array}\right.
$$

\noindent Thus, two paths in $G(X,R)$ identified by $\equi{(X,R)}$ are sent by $\phi$ on two paths in $G(\Sigma^3(X,R))$ identified by $\equi{\Sigma^3(X,R)}$. The inclusion of $\psi(\equi{\Sigma^3(X,R)})$ into $\equi{(X,R)}$ is proved similarly, starting from the last two equations, in the case $c=0$, and moving upwards to a square whose paths are identified by $\equi{(X,R)}$.

Now, let us fix a $3$-polygraph $\Sigma=(\star,\emptyset,\Sigma_2,\Sigma_3)$. The equality $\Sigma^3(\Nr(\Sigma))=\Sigma$ is once again immediate. Since $\Nr(\Sigma)$ is a Petri net, we know that $G(\Nr(\Sigma))$ is isomorphic to $G(\Sigma^3(\Nr(\Sigma)))$, which is the same as $G(\Sigma)$. Furthermore, this graph isomorphism is defined the same way as $\phi$ and $\psi$ in the first part of the proof. Hence $\phi(\equi{\Nr(\Sigma)})$ is equal to $\equi{\Sigma}$. We apply $\psi$ to get the equality between $\equi{\Nr(\Sigma)}$ and $\psi(\equi{\Sigma})$. 

\findem\end{dem}

\noindent This result allows the informal statement "Petri nets \emph{are exactly} the $3$-polygraphs with one $0$-cell and no $1$-cell" for the following reasons:
\begin{enumerate}
\item[-] There is a correspondance between the presentations, given by the interpretation of places as $2$-cells and of transitions as $3$-cells.
\item[-] Both presentations generate the same reduction graph, so that each one can simulate the evolutions of the other one.
\item[-] There is a correspondance between the congruences that identify, in each graph, the paths that only differ by the order of application of the same transitions/3-cells. 
\end{enumerate}

\noindent Another, more categorical way to formulate this correspondance is to say that the category $\Gr(X,R)$ generated by a Petri net is isomorphic to the category whose objects and arrows are respectively the $2$-arrows and $3$-arrows of $\mon{\Sigma^3(X,R)}$.

\section*{Comments and future directions}
\emptysectionmark{Comments and future directions}

We have proved that Petri nets have two natural interpretations in terms of polygraphs. Let us informally compare them.

The first one, using a $2$-polygraph, is really convenient to use, since the circuit-like representation is now well-understood and user-friendly. The only difficulty comes with the explicit permutations: one has to choose a way to identify two paths that only differ by permutations. We have discussed possible starting points in order to reach a solution for this issue. And, as we have seen, this is non trivial and is postponed to further work. Nonetheless, this is an important new challenge for $3$-dimensional rewriting, since the polygraphs involved provide a new class of rather different examples.

The second polygraphic interpretation we have studied, using a $3$-dimensional polygraph, provides, at least theoretically, a better description of the intrinsic algebraic structure of Petri nets: they do not require any extra cell, apart from the ones given with the Petri nets. However, these objects are hard to handle for the moment and this mainly comes from the lack of graphical representations: indeed, the first ones have been constructed in [Guiraud 2005] to represent classical proofs, but they remain hard to produce and handle in a convenient way. For that reason, part of the future work will concern these $3$-dimensional representations: the goals are to improve the ones already known, to automatize their production and, maybe, to search for other ones. In the case of Petri nets, the representations should be really interesting since their shape will strangely be close to diagrams used in superstring theory to represent interactions between superstrings.

Let us finish by a more general comment on polygraphs. The results presented here constitute another clue of the expressive power of polygraphs in theoretical computer science, proof theory and universal algebra. Indeed, it is already known that polygraphs generalize word and term rewriting systems, equational presentations of algebraic structures, Reidemeister moves on knots and tangles, formal proofs of classical logic. The interested reader can find more information about the translations of all these objects into polygraphs in the following documents: [Burroni 1993], [Lafont 2003], [Métayer 2003], [Guiraud 2004(T), 2004, 2005]. 

\vfill
\begin{flushright}
\begin{minipage}{100mm}
\emph{I wish to thank Albert Burroni and Yves Lafont for many discussions and advices and the referees for their comments that have helped to improve this document.}
\end{minipage}
\end{flushright}

\section*{References}
\emptysectionmark{References}

\noindent\textsc{Franz Baader}, \textsc{Tobias Nipkow} \\
\indent\emph{Term rewriting and all that}, Cambridge University Press, 1998.  

\bigskip
\noindent\textsc{John Carlos Baez}, \textsc{James Dolan} \\
\indent\emph{Categorification}, ArXiv preprint, 1998.  

\bigskip
\noindent\textsc{Albert Burroni} \\ 
\indent\emph{Higher-dimensional word problems with applications to equational logic}, \\
\indent\indent Theoretical Computer Science 115(1), 1993. 

\bigskip
\noindent\textsc{Olga Caprotti, Alois Ferscha, Hoon Hong} \\
\indent\emph{Reachability test in Petri nets by Gröbner bases}, RISC report series 95(03), 1995. 

\bigskip
\noindent\textsc{Angie Chandler, Anne Heyworth} \\
\indent\emph{Gröbner bases as a tool for Petri net analysis}, Proceedings SCI 2001. 

\bigskip
\noindent\textsc{Yves Guiraud} \\
\indent\emph{Présentations d'opérades et systèmes de réécriture}, Thèse de doctorat, 2004(T). \\ 
\indent\emph{Termination orders for $3$-dimensional rewriting}, \\
\indent\indent To appear in Journal of Pure and Applied Algebra (2004). \\
\indent\emph{The three dimensions of proofs}, To appear in Annals of Pure and Applied Logic (2005). 

\bigskip
\noindent\textsc{Yves Lafont} \\
\indent\emph{Towards an algebraic theory of boolean circuits}, Journal of Pure and Applied Algebra 184, 2003. 

\bigskip
\noindent\textsc{Saunders MacLane} \\
\indent\emph{Categories for the working mathematician}, Springer, second edition 1998. 

\bigskip
\noindent\textsc{François Métayer} \\
\indent\emph{Resolutions by polygraphs}, Theory and Applications of Categories 11(7), 2003. 

\bigskip
\noindent\textsc{Tadao Murata} \\
\indent\emph{Petri nets: properties, analysis and applications}, Proceedings IEEE 77(4), 1989.

\end{document}